\documentclass[12pt]{article}
\usepackage{amssymb,amsmath,amsthm}
\usepackage[all]{xy}           

\makeatletter
\begin{document}

\newtheorem{theorem}{Theorem}[section]
\newtheorem{prop}[theorem]{Proposition}
\newtheorem{defn}[theorem]{Definition}
\newtheorem{lemma}[theorem]{Lemma}
\newtheorem{coro}[theorem]{Corollary}
\newtheorem{prop-def}{Proposition-Definition}[section]
\newtheorem{claim}{Claim}[section]
\newtheorem{propprop}{Proposed Proposition}[section]
\newtheorem{conjecture}{Conjecture}
\newcommand{\nc}{\newcommand}

\nc{\bin}[2]{ (_{\stackrel{\scs{#1}}{\scs{#2}}})}  
\nc{\binc}[2]{(\!\! \begin{array}{c} \scs{#1}\\
    \scs{#2} \end{array}\!\!)}  
\nc{\bincc}[2]{  ( {\scs{#1} \atop
    \vspace{-1cm}\scs{#2}} )}  
\nc{\bs}{\bar{S}}
\nc{\la}{\longrightarrow}
\nc{\rar}{\rightarrow}
\nc{\dar}{\downarrow}
\nc{\dap}[1]{\downarrow \rlap{$\scriptstyle{#1}$}}
\nc{\defeq}{\stackrel{\rm def}{=}}
\nc{\dis}[1]{\displaystyle{#1}}
\nc{\dotcup}{\ \displaystyle{\bigcup^\bullet}\ }
\nc{\hcm}{\ \hat{,}\ }
\nc{\hts}{\hat{\otimes}}
\nc{\hcirc}{\hat{\circ}}
\nc{\lleft}{[}
\nc{\lright}{]}
\nc{\curlyl}{\left \{ \begin{array}{c} {} \\ {} \end{array}
    \right .  \!\!\!\!\!\!\!} 
\nc{\curlyr}{ \!\!\!\!\!\!\!
    \left . \begin{array}{c} {} \\ {} \end{array}
    \right \} }
\nc{\longmid}{\left | \begin{array}{c} {} \\ {} \end{array}
    \right . \!\!\!\!\!\!\!}
\nc{\ora}[1]{\stackrel{#1}{\rar}}
\nc{\ola}[1]{\stackrel{#1}{\la}}
\nc{\scs}[1]{\scriptstyle{#1}}
\nc{\mrm}[1]{{\rm #1}}
\nc{\dirlim}{\displaystyle{\lim_{\longrightarrow}}\,}
\nc{\invlim}{\displaystyle{\lim_{\longleftarrow}}\,}
\nc{\mvp}{\vspace{0.3cm}}
\nc{\tk}{^{(k)}}
\nc{\tp}{^\prime}
\nc{\ttp}{^{\prime\prime}}
\nc{\svp}{\vspace{2cm}}
\nc{\vp}{\vspace{8cm}}
\nc{\proofend}{$\blacksquare$ \vspace{0.3cm}}
\nc{\modg}[1]{\!<\!\!{#1}\!\!>}
\nc{\intg}[1]{F_C(#1)}
\nc{\lmodg}{\!<\!\!}
\nc{\rmodg}{\!\!>\!}
\nc{\cpi}{\widehat{\Pi}}
\nc{\sha}{{\mbox{\cyr X}}}  
\nc{\shpr}{\diamond}    
\nc{\shprc}{\shpr_c}
\nc{\labs}{\mid\!}
\nc{\rabs}{\!\mid}

\nc{\ann}{\mrm{ann}}
\nc{\Aut}{\mrm{Aut}}
\nc{\can}{\mrm{can}}
\nc{\colim}{\mrm{colim}}
\nc{\Cont}{\mrm{Cont}}
\nc{\rchar}{\mrm{char}}
\nc{\cok}{\mrm{coker}}
\nc{\dtf}{{R-{\rm tf}}}
\nc{\dtor}{{R-{\rm tor}}}
\renewcommand{\det}{\mrm{det}}
\nc{\Div}{{\mrm Div}}
\nc{\End}{\mrm{End}}
\nc{\Ext}{\mrm{Ext}}
\nc{\Fil}{\mrm{Fil}}
\nc{\Frob}{\mrm{Frob}}
\nc{\Gal}{\mrm{Gal}}
\nc{\GL}{\mrm{GL}}
\nc{\Hom}{\mrm{Hom}}
\nc{\hsr}{\mrm{H}}
\nc{\hpol}{\mrm{HP}}
\nc{\id}{\mrm{id}}
\nc{\im}{\mrm{im}}
\nc{\incl}{\mrm{incl}}
\nc{\length}{\mrm{length}}
\nc{\mchar}{\rm char}
\nc{\mpart}{\mrm{part}}
\nc{\ql}{{\QQ_\ell}}
\nc{\qp}{{\QQ_p}}
\nc{\rank}{\mrm{rank}}
\nc{\rcot}{\mrm{cot}}
\nc{\rdef}{\mrm{def}}
\nc{\rdiv}{{\rm div}}
\nc{\rtf}{{\rm tf}}
\nc{\rtor}{{\rm tor}}
\nc{\res}{\mrm{res}}
\nc{\SL}{\mrm{SL}}
\nc{\Spec}{\mrm{Spec}}
\nc{\tor}{\mrm{tor}}
\nc{\Tr}{\mrm{Tr}}
\nc{\tr}{\mrm{tr}}

\nc{\ab}{\mathbf{Ab}}
\nc{\Alg}{\mathbf{Alg}}
\nc{\Bax}{\mathbf{Bax}}
\nc{\bfk}{{\bf k}}
\nc{\bfone}{{\bf 1}}
\nc{\changes}{\marginpar{\bf Changes!}}
\nc{\detail}{\marginpar{\bf More detail}
    \noindent{\bf Need more detail!}
    \svp}
\nc{\Diff}{\mathbf{Diff}}   
\nc{\gap}{\marginpar{\bf Incomplete}\noindent{\bf Incomplete!!}
    \svp}
\nc{\FMod}{\mathbf{FMod}}
\nc{\Int}{\mathbf{Int}}
\nc{\Mon}{\mathbf{Mon}}
\nc{\remark}{\noindent{\bf Remark: }}
\nc{\remarks}{\noindent{\bf Remarks: }}
\nc{\Rep}{\mathbf{Rep}}
\nc{\Rings}{\mathbf{Rings}}
\nc{\Sets}{\mathbf{Sets}}
\nc{\bill}[1]{\marginpar{\bf To Bill}\noindent{\bf To Bill:}
    {\tt #1}\\ }
\nc{\li}[1]{\marginpar{\bf To Li}\noindent{\bf To Li:}
    {\tt #1}\\ }

\nc{\BA}{{\Bbb A}}
\nc{\CC}{{\Bbb C}}
\nc{\DD}{{\Bbb D}}
\nc{\EE}{{\Bbb E}}
\nc{\FF}{{\Bbb F}}
\nc{\GG}{{\Bbb G}}
\nc{\HH}{{\Bbb H}}
\nc{\LL}{{\Bbb L}}
\nc{\NN}{{\Bbb N}}
\nc{\QQ}{{\Bbb Q}}
\nc{\RR}{{\Bbb R}}
\nc{\TT}{{\Bbb T}}
\nc{\VV}{{\Bbb V}}
\nc{\ZZ}{{\Bbb Z}}


\nc{\cala}{{\cal A}}
\nc{\calc}{{\cal C}}
\nc{\cald}{\mathcal{D}}
\nc{\cale}{{\cal E}}
\nc{\calf}{{\cal F}}
\nc{\calg}{{\cal G}}
\nc{\calh}{{\cal H}}
\nc{\cali}{{\cal I}}
\nc{\call}{{\cal L}}
\nc{\calm}{{\cal M}}
\nc{\caln}{{\cal N}}
\nc{\calo}{{\cal O}}
\nc{\calp}{{\cal P}}
\nc{\calr}{{\cal R}}
\nc{\calt}{{\cal T}}
\nc{\calw}{{\cal W}}
\nc{\calx}{{\cal X}}
\nc{\CA}{\mathcal{A}}

\nc{\fraka}{{\frak a}}
\nc{\frakB}{{\frak B}}
\nc{\frakm}{{\frak m}}
\nc{\frakp}{{\frak p}}

\font\cyr=wncyr10

\title{ 
Baxter Algebras and Shuffle Products
\thanks{The first author is supported in part by NSF grant
    \#DMS 97-96122. 
    MSC Numbers: Primary 16A06, 47B99. 
    Secondary 13A99,16W99.}}
\author{Li Guo and 
William Keigher
\\
Department of Mathematics and Computer Science\\
Rutgers University\\
Newark, NJ 07102 \\
(liguo@newark.rutgers.edu)\\
(keigher@newark.rutgers.edu)
}
\date{}
\maketitle

\setcounter{section}{0}

\section{Introduction}

In this paper we generalize the well-known construction
of shuffle product algebras by using mixable shuffles,
and prove that any free Baxter algebra is isomorphic
to a mixable shuffle product algebra.
This gives an explicit construction of the free Baxter
algebra, extending the work of Rota~\cite{Rot1}
and Cartier~\cite{Ca}.

In an important paper published in 1958, 
Ree~\cite{Re} constructed algebras in which the product is
expressed in terms of shuffles. 
He was motivated by Chen's work on iterated integrals of
paths~\cite{Ch3}, where this shuffle product 
is derived from the integration by parts formula
\begin{equation}
 \int _0^x f(t)dt \int_0^x g(t)dt
= \int_0^x f(t)( \int_0^t g(s)ds) dt
+ \int_0^x g(t)( \int_0^t f(s)ds) dt.
\label{eq:int}
\end{equation}
Subsequently, shuffle product constructions have been
studied extensively and have found 
applications in many areas of pure and applied
mathematics.

In another important paper published in 1960,
Baxter~\cite{Ba} considered operators $P$ that satisfies the
identity
\[ P(x)P(y)+P(xy)=P(xP(y))+P(yP(x))\]
and used this identity to study the theory of fluctuations.
Rota studied Baxter's operators from
an algebraic point of view and defined a Baxter algebra
to be an algebra $A$ with an operator $P$ satisfying the identity
\[ P(x)P(y)+q P(xy)=P(xP(y))+P(yP(x))\]
for some fixed $q$ in the base ring of $A$. 
Rota~\cite{Rot1} and Cartier~\cite{Ca} gave explicit constructions
of the free Baxter algebra on a set $X$ in the case
when $q=1$.

It is easy to see that
the product in Baxter algebras when $q=0$ 
can be described by shuffle products, and
the shuffle product algebras considered by Chen and Ree 
have a canonical Baxter operator.
However, there does not appear to be an explicit and systematic
study of this connection between
Baxter algebras and shuffle products except a remark
in a recent paper of Rota~\cite{Rot4} implying
such a connection.

The motivation for the current paper came from the desire of
developing a theory that
is ``dual" to the beautiful theory of differential algebras
obtained by Ritt~\cite{Ri} and Kolchin~\cite{Ko}.
Since a differential algebra is an algebra with an operator
that satisfies the Leibniz product rule, it seems
natural to study ``integration algebras",
i.e., algebras with an operator that satisfies
an identity similar to the one in equation~(\ref{eq:int}). 
Of course, these are just Baxter algebras where $q=0$.
Unaware of the above mentioned work on Baxter algebras,
we gave a description of the
free integration algebra by using shuffle products. 

It was Rota who pointed out to us the earlier work
on free Baxter algebras and 
suggested that we extend our shuffle product description
of free integration algebras to Baxter algebras.
This is carried out in this paper, 
by making use of a modified shuffle product,
called the mixable shuffle product.
Thus we not only consider shuffles of two vectors,
but also shuffles in which certain
components of a shuffle will ``merge". 
This enables us to construct the mixable shuffle
product algebras, generalizing the classical construction
of shuffle product algebras,
and to give a more intuitive
and constructive description of the free Baxter algebra. 
Also, our description is of
free Baxter algebras on any commutative algebra
with any value $q$.
Furthermore, 
the free Baxter algebra of Cartier or Rota
is in the category of Baxter algebras not necessarily having
an identity, while the free Baxter algebra we consider
is in the category
of Baxter algebras with an identity. 
When specialized to the case considered by Rota or Cartier, 
the free Baxter algebra we construct
contains their free Baxter algebra as a
sub-Baxter algebra in the category of Baxter algebras not
necessarily having an identity. 

The mixable shuffle product of Baxter algebras
allows us to study in more detail the properties of
Baxter algebras. This will be the subject of a forthcoming
paper. 
Shuffle products have occurred in many other fields and contexts,
such as Hopf algebras, algebraic $K$-theory, algebraic
topology and combinatorics, as well as in computational
mathematics and applied mathematics; 
see, for example, ~\cite{Ch1,Ch2,CL,Koc, Re,Ro,Sw}.
One might hope that the mixable shuffle introduced here will
provide interesting and useful generalizations to those
theories.
One might also hope that the mixable shuffle product 
will be useful for the applications of Baxter algebras,
such as those considered by Baxter and Rota.

In this paper, we first give a brief summary of
basic definitions and properties of Baxter algebras
in section~\ref{sec:def}. 
We then make a careful study of mixable shuffles
in section~\ref{sec:shuf}.
In section~\ref{sec:free} we apply properties of mixable
shuffles to construct mixable shuffle product algebras
and to prove that any free Baxter algebra is isomorphic
to such an algebra. 
In this section, we also briefly consider variations of
the free construction by examining the free
Baxter algebra on a set, on a commutative monoid
and on a module.
The relation between the free Baxter algebras we construct 
and the free Baxter algebras constructed by Rota and Cartier
is described in section~\ref{sec:RC}. 
We conclude by considering the special case of
the free Baxter algebra on the empty set.

\section{Definitions and basic properties}
\label{sec:def}
In this paper, any ring $R$ is commutative with
identity element $\bfone_R$.
All notation will be standard unless otherwise
noted.  In particular, we write $\NN$ for the additive monoid of
natural numbers $\{0,1,2,\ldots\}$ and 
$\NN_+=\{ n\in \NN\mid n>0\}$ for the positive integers. 
Also, $\bincc nk$ will denote the
usual binomial coefficient defined for any $n, k \in \NN$ with $k
\leq n$ by  $\bincc nk = n!/(k!(n - k)!)$.

\begin{defn}
Let $C$ be a ring, $q\in C$,
and let $R$ be a $C$-algebra.
A {\bf Baxter operator on $R$ over $C$} is a $C$-module
endomorphism $P$ of $R$ satisfying 
\begin{equation}
 P(x)P(y)+qP(xy)=P(xP(y))+P(yP(x)),\ x,\ y\in R.
\label{eq:bax1}
\end{equation}
\end{defn}
We will find it convenient for the remainder of this
paper to write equation~(\ref{eq:bax1}) in the form
\begin{equation}
 P(x)P(y)=P(xP(y))+P(yP(x))+\lambda P(xy),
\label{eq:bax2}
\end{equation}
so that our $\lambda$ is $-q$.
We will also say that $P$ has weight $\lambda$.

\begin{defn}
    A {\bf Baxter C-algebra of weight $\lambda$}
is a pair $(R,P)$ where $R$ is a $C$-algebra
and $P$ is a Baxter operator of weight $\lambda$
on $R$ over $C$.
\end{defn}
If the meaning of $\lambda$ is clear, we will
suppress $\lambda$ from the notation.

Note that the mapping $0 : R \rar R$ defined by $0(r) = 0$ for all 
$r \in R$ is trivially a Baxter operator
on $R$ over $R$, for any ring $R$.
Hence every $C$-algebra can be viewed as
a Baxter $C$-algebra.
Definitions of basic concepts for $C$-algebras can 
be similarly defined for Baxter $C$-algebras.
In particular,
let $(R,P)$ and $(S,Q)$ be two Baxter $C$-algebras
of weight $\lambda$.
A {\bf homomorphism of Baxter $C$-algebras}
$f:(R,P)\rar (S,Q)$ is a homomorphism
$f: R \rar S$ of $C$-algebras with the property that
$ f(P(x))=Q(f(x))$
for all $x \in R$.

Let $\Bax_{C,\lambda}$ denote the category of Baxter
$C$-algebras of weight $\lambda$. 
A {\bf Baxter ideal} of $(R,P)$ is an ideal $I$ of $R$ such
that $P(I)\subseteq I$.
For a Baxter ideal $I$ of $(R,P)$,
the {\bf quotient Baxter $C$-algebra} is
the quotient algebra $R/I$, together with the
$C$-linear endomorphism
$\bar{P}:R/I\rar R/I$ induced from $P$. 
If $f:(R,P)\rar (S,Q)$ is in $\Bax_{C,\lambda}$, then
$\ker f$ is a Baxter ideal of $R$,
and $\im f$, with the restriction of $Q$, is a Baxter
sub-$C$-algebra of $(S,Q)$.

\section{Mixable shuffle products}
\label{sec:shuf}
This section is a preparation for the next section,
where we will give a description of
free Baxter algebras in terms of mixable shuffles. 
We start with a study of mixable shuffles in the context of
permutations.
We then apply this study to mixable shuffles of vectors.
Finally we consider mixable shuffles of tensors,
using the mixable shuffle of vectors
as a ``generic" form. 
For the purpose of providing a solid foundation for later
applications,
we give full details of the proofs, even though some of
them might be intuitively clear to the expert.
So some readers might just want to look at the definitions
and results and move on to the next section.

\subsection{Permutation shuffles}
For $m,n\in \NN_+$, 
define the set of {\bf $(m,n)$-shuffles} by 
\begin{eqnarray*}
\lefteqn{ S(m,n)=} \\
&& \left \{ \sigma\in S_{m+n}
    \begin{array}{ll} {} \\ {} \end{array} \right .
\left | 
\begin{array}{l}
\sigma^{-1}(1)<\sigma^{-1}(2)<\ldots<\sigma^{-1}(m),\\
\sigma^{-1}(m+1)<\sigma^{-1}(m+2)<\ldots<\sigma^{-1}(m+n)
\end{array}
\right \}.
\end{eqnarray*}
Equivalently, 
\begin{eqnarray*}
\lefteqn{ S(m,n)=} \\
 && \left \{ \sigma\in S_{m+n}
    \begin{array}{l} {} \\ {} \end{array} \right .
\left | \begin{array}{ll}
{\rm \ if\ }
    1\leq \sigma(r)<\sigma(s)\leq m \\
    {\rm \ or\ } m+1\leq \sigma(r)<\sigma(s)\leq m+n 
    {\rm \ then\ } r<s
    \end{array}
\right \}.
\end{eqnarray*}

A pair of indices $(k, k+1)$,\ $1\leq k< m+n$ is 
called an {\bf admissible pair} for
an $(m,n)$-shuffle $\sigma\in S(m,n)$
if $\sigma(k)\leq m<\sigma(k+1)$. 
Denote $\calt^\sigma$ for the set of admissible pairs for $\sigma$. 
For $\sigma\in S(m,n)$ and $T\subseteq \calt^\sigma$, call the pair
$(\sigma,T)$ a {\bf mixable $(m,n)$-shuffle}.
When $T=\phi$, $(\sigma,T)$ is identified with
the shuffle $\sigma$. 
Denote $\bs (m,n)$ for the set of mixable $(m,n)$-shuffles.
So
\[ \bs (m,n)=\{ (\sigma,T)\mid \sigma\in S(m,n),\
    T\subseteq \calt^\sigma\}. \]
Also denote
\[ s(m,n)=\mid \bs(m,n)\mid .\]

Intuitively, a $(m,n)$-shuffle is a permutation $\sigma$
of $\{1,\ldots,m, m+1,\ldots, m+n\}$ such that
the order of $\{1,\ldots,m\}$ and $\{m+1,\ldots,m+n\}$
are preserved in
$\{ \sigma(1),\ldots,\sigma(m),\sigma(m+1),\ldots,\sigma(m+n)\}$.
Further, a mixable $(m,n)$-shuffle $(\sigma,T)$
is a $(m,n)$-shuffle $\sigma$ in which pairs of indices
from $T$ represent positions where ``merging" will occur.
The precise meaning of ``merging" will be described
later.

For example, 
$\sigma=\left ( \begin{array}{ccc} 1 & 2& 3 \\
    1&3&2 \end{array} \right ) $
is a $(2,1)$-shuffle.
The pair $(1,2)$ is an admissible pair for $\sigma$.
Note that
$\tau=\left ( \begin{array}{ccc} 1 & 2& 3 \\
    3&1& 2 \end{array} \right ) $
is also a $(2,1)$-shuffle, but has no admissible pairs.

Denote
\begin{eqnarray*}
 \bs_{1,0}(m,n)&=& \{(\sigma,T)\in \bs(m,n) \mid (1,2)\not\in T,
    \sigma^{-1}(1)=1 \},\\
 \bs_{0,1}(m,n)&=& \{(\sigma,T)\in \bs(m,n)\mid (1,2)\not\in T,
    \sigma^{-1}(m+1)=1 \},\\
 \bs_{1,1}(m,n)&=&\{ (\sigma,T)\in \bs(m,n)\mid (1,2)\in T\}.
\end{eqnarray*}
Then we clearly have
\begin{equation} \bs (m,n) = \bs_{1,0}(m,n) \dotcup
    \bs_{0,1}(m,n) \dotcup \bs_{1,1}(m,n)
\label{eq:dec}
\end{equation}
since, by the definition of $S(m,n)$, one
and only one of $\sigma^{-1}(1)$ and
$\sigma^{-1}(m+1)$ equals 1. 
Similarly we will define the set of {\bf $(m,n,\ell)$-shuffles} by 

\begin{eqnarray*}
\lefteqn{ S(m,n,\ell)=} \\
 && \left \{ \sigma\in S_{m+n+\ell}
    \begin{array}{ll} {} \\ {} \end{array} \right .
\left | 
\begin{array}{l}
\sigma^{-1}(1)<\sigma^{-1}(2)<\ldots<\sigma^{-1}(m),\\
\sigma^{-1}(m+1)<\sigma^{-1}(m+2)<\ldots<\sigma^{-1}(m+n), \\
\sigma^{-1}(m+n+1)<\sigma^{-1}(m+n+2) \\
\ \ \ \ <\ldots<\sigma^{-1}(m+n+\ell)
\end{array}
\right \}.
\end{eqnarray*}
For $\sigma\in S(m,n,\ell)$,  define 
\[  \calt^\sigma  = \calt^\sigma_{1,1,0} \dotcup
    \calt^\sigma_{1,0,1} \dotcup
    \calt^\sigma_{0,1,1} \dotcup
    \calt^\sigma_{1,1,1},\]
where 
\begin{eqnarray*}
\calt^\sigma_{1,1,0}\!\!\! &= &\!\!\!  \curlyl (k,k+1) \longmid
    \begin{array}{l}
    1\leq k <m+n+\ell-1, 
     \sigma(k)\leq m<\sigma (k+1)\\
    {\rm and\ either\ } 
    \sigma (k+1)> m+n
    {\rm\ or\ }
    m+n\geq \sigma(k+2)
    \end{array}
    \curlyr,  \\
\calt^\sigma_{1,0,1} &= &  \curlyl (k,k+1) \longmid
    \begin{array}{l}
    1\leq k <m+n+\ell, 
     \sigma(k)\leq m,\
    m+n< \sigma (k+1)
    \end{array}
    \curlyr,  \\
\calt^\sigma_{0,1,1} &= &  \curlyl (k,k+1) \longmid
    \begin{array}{l}
    1< k <m+n+\ell, 
     \sigma(k) \leq m+n< \sigma(k+1)\\
     {\rm\  and\ either\ }
     \sigma(k-1)>m
     {\rm\ or\ }
     m\geq \sigma(k)
    \end{array}
    \curlyr,  \\
\calt^\sigma_{1,1,1} &= & 
    \curlyl (k,k+1,k+2)\longmid
    \begin{array}{l} 1\leq k <m+n+\ell-1,\\
    \sigma(k) \leq m <\sigma(k+1) \\
    \ \ \ \leq m+n <\sigma(k+2)
    \end{array} \curlyr. 
\end{eqnarray*}
For $(\sigma,T)\in \bs(m,n,\ell)$, define
\begin{equation}
\deg T = \mid T\cap \calt^\sigma_{1,1,0} \mid +
    \mid T\cap \calt^\sigma_{1,0,1} \mid +
    \mid T\cap \calt^\sigma_{0,1,1} \mid +
    2\mid T\cap \calt^\sigma_{1,1,1} \mid.
\label{eq:deg3}
\end{equation}
Denote
\[ \bs (m,n,\ell) =
    \{ (\sigma,T)\mid \sigma\in S(m,n,\ell),
    T\in \calt_\sigma\}, \]
and
\[ s(m,n,\ell) = \mid \bs(m,n,\ell)\mid.\]

\begin{prop}
\label{prop:num}
Let $m,n\in \NN_+$. 
\begin{enumerate}
\item
$ s(m,n)=s(m-1,n)+ s(m,n-1)+ s(m-1,n-1).$
\item
$ \mid \{ (\sigma,T)\in \bs(m,n) \mid\   \mid T\mid =i\} \mid
 = \bincc{m+n-i}{n}\bincc{n}{i}.$
\item
$ s(m,n) =
    \displaystyle{\sum_{i=0}^n \bincc{m+n-i}{n} \bincc{n}{i}}.$
\item   
$ s(m,n,\ell)=
    \displaystyle{ \sum_{k=0}^{n+\ell}
    \sum_{i=0}^n
    \bincc{m+n+\ell-k}{\ell} \bincc{\ell}{k-i}
    \bincc{m+n-i}{n} \bincc{n}{i}.} $
\end{enumerate} 
\end{prop}  

\proof
1.
This follows from equation~(\ref{eq:dec}).

2.
Denote the set on the left hand side of the
equation by $\bs^{(i)}(m,n)$. 
We use induction on $m+n$, with $m,\ n\geq 1$.
When $m+n=2$, this can be verified directly.
In general,
define
$\bs^{(i)}_{(1,0)}(m,n) = \bs^{(i)}(m,n)\cap
    \bs_{(1,0)}(m,n)$
and similarly for
$\bs^{(i)}_{(0,1)}(m,n)$
and $\bs^{(i)}_{(1,1)}(m,n)$. From equation~(\ref{eq:dec}),
we obtain,
\begin{eqnarray*}
 \bs^{(i)} (m,n)
&=& \bs^{(i)}_{1,0}(m,n) \dotcup
    \bs^{(i)}_{0,1}(m,n) \dotcup \bs^{(i)}_{1,1}(m,n)\\
&\cong & \bs^{(i)}(m-1,n) \dotcup \bs^{(i)}(m,n-1)
    \dotcup \bs^{(i-1)} (m-1,n-1). 
\end{eqnarray*}
Here $A\cong B$ means that the two sets $A$ and $B$ have
the same cardinality. 
By the inductive assumption,
the size of the right hand side is
\[ \bincc{m+n-1-i}{m-1}\bincc{m-1}{i}
+ \bincc{m+n-1-i}{m}\bincc{m}{i}
+\bincc{m+n-2-i}{m-1}\bincc{m-1}{i-1}.
\]
Applying Pascal's identity, we see that
the last sum is also the value of
$\bincc{m+n-i}{m}\bincc{m}{i}$.

3.
This follows from part 2 by summing over $i$
for $i=0,\ldots,n$.

4.
The proof is similar to the proof of part 3. 
For $0\leq k\leq n+\ell$, denote
\[ \bs^{(k)}(m,n,\ell) =
    \curlyl (\sigma,T)\in \bs(m,n,\ell) \mid
    \deg\  T=k \curlyr.\]
We will use induction on $j=m+n+\ell,\  m,\ n,\ \ell\geq 1$ to prove

\begin{equation}
 \mid \bs^{(k)}(m,n,\ell)\mid
 = \displaystyle{
    \sum_{i=0}^n }\bincc{m+n+\ell-k}{\ell}\bincc{\ell}{k-i}
    \bincc{m+n-i}{n}\bincc{n}{i}.
\label{eq:3sh}
\end{equation}

When $m+n+\ell=3$, the equation can be verified directly.
Now assume that the equation holds for $m+n+\ell<j$
and consider the case when $m+n+\ell=j$.
Denote 
\begin{eqnarray*}
 \bs_{1,0,0}(m,n,\ell)&=& \{(\sigma,T)\in \bs(m,n) \mid
    \sigma^{-1}(1)=1, (1,2)\not\in T \},\\
 \bs_{0,1,0}(m,n,\ell)&=& \{(\sigma,T)\in \bs(m,n) \mid 
    \sigma^{-1}(m+1)=1,(1,2)\not\in T \},\\
 \bs_{0,0,1}(m,n,\ell)&=& \{(\sigma,T)\in \bs(m,n) \mid 
    \sigma^{-1}(m+n+1)=1, (1,2)\not\in T \},\\
 \bs_{1,1,0}(m,n,\ell)&=& \{(\sigma,T)\in \bs(m,n) \mid
    (1,2)\in T\cap \calt^\sigma_{1,1,0},
    (1,2,3)\not\in T \},\\
 \bs_{1,0,1}(m,n,\ell)&=& \{(\sigma,T)\in \bs(m,n) \mid
    (1,2)\in T\cap \calt^\sigma_{1,0,1},
    (1,2,3)\not\in T \},\\
 \bs_{0,1,1}(m,n,\ell)&=& \{(\sigma,T)\in \bs(m,n) \mid
    (1,2)\in T\cap \calt^\sigma_{0,1,1},
    (1,2,3)\not\in T \},\\
 \bs_{1,1,1}(m,n,\ell)&=& \{(\sigma,T)\in \bs(m,n) \mid
    (1,2,3)\in T \}. 
\end{eqnarray*}
Also denote
\[ \bs^{(k)}_{u,v,w}(m,n,\ell)
    =\bs_{u,v,w}(m,n,\ell) \cap \bs^{(k)}(m,n,\ell)\]
for $u,v,w=0$ or 1. 
It follows from the definition that 
\begin{eqnarray*}
\lefteqn{ \bs (m,n,\ell)=}\\
&&  \bs_{1,0,0}(m,n,\ell) \dotcup
    \bs_{0,1,0}(m,n,\ell) \dotcup
    \bs_{0,0,1}(m,n,\ell) \\
&&  \dotcup
    \bs_{1,1,0}(m,n,\ell) \dotcup
    \bs_{1,0,1}(m,n,\ell) \dotcup
    \bs_{0,1,1}(m,n,\ell) \dotcup 
    \bs_{1,1,1}(m,n,\ell).
\end{eqnarray*}
So
\begin{eqnarray*}
\lefteqn{ \bs^{(k)} (m,n,\ell)=}\\
&&  \bs^{(k)}_{1,0,0}(m,n,\ell) \dotcup
    \bs^{(k)}_{0,1,0}(m,n,\ell) \dotcup
    \bs^{(k)}_{0,0,1}(m,n,\ell) \\
&&  \dotcup 
    \bs^{(k)}_{1,1,0}(m,n,\ell) \dotcup
    \bs^{(k)}_{1,0,1}(m,n,\ell) \dotcup
    \bs^{(k)}_{0,1,1}(m,n,\ell) \dotcup 
    \bs^{(k)}_{1,1,1}(m,n,\ell)\\
&\cong &
    \bs^{(k)}_{1,0,0}(m-1,n,\ell) \dotcup
    \bs^{(k)}_{0,1,0}(m,n-1,\ell) \dotcup
    \bs^{(k)}_{0,0,1}(m,n,\ell-1) \\
&&  \dotcup
    \bs^{(k-1)}_{1,1,0}(m-1,n-1,\ell) \dotcup
    \bs^{(k-1)}_{1,0,1}(m-1,n,\ell-1) \\
&&  \dotcup
    \bs^{(k-1)}_{0,1,1}(m,n-1,\ell-1) 
    \dotcup
    \bs^{(k-2)}_{1,1,1}(m-1,n-1,\ell-1). 
\end{eqnarray*}
Here again $\cong$ stands for a bijection between sets.
Taking cardinalities and applying the induction hypothesis,
we see that the left hand side of equation~(\ref{eq:3sh}) is 
\allowdisplaybreaks{
\begin{eqnarray*}
& & \displaystyle{
    \sum_{i=0}^n }\bincc{m+n+\ell-k-1}{\ell}\bincc{\ell}{k-i}
    \bincc{m+n-i-1}{n}\bincc{n}{i} \\
& + &\displaystyle{
    \sum_{i=0}^{n-1} }\bincc{m+n+\ell-k-1}{\ell}\bincc{\ell}{k-i}
    \bincc{m+n-i-1}{n-1}\bincc{n-1}{i} \\
& +&  \displaystyle{
    \sum_{i=0}^n }\bincc{m+n+\ell-k-1}{\ell-1}\bincc{\ell-1}{k-i}
    \bincc{m+n-i-1}{n}\bincc{n}{i} \\
& +&  \displaystyle{
    \sum_{i=0}^{n-1} }\bincc{m+n+\ell-k-1}{\ell}
        \bincc{\ell}{k-i-1}
    \bincc{m+n-i-2}{n-1}\bincc{n-1}{i} \\
& +&  \displaystyle{
    \sum_{i=0}^n }\bincc{m+n+\ell-k-1}{\ell-1}
        \bincc{\ell-1}{k-i-1}
    \bincc{m+n-i-1}{n}\bincc{n}{i} \\
& +&  \displaystyle{
    \sum_{i=0}^{n-1} }\bincc{m+n+\ell-k-1}{\ell}
        \bincc{\ell-1}{k-i-1}
    \bincc{m+n-i-1}{n-1}\bincc{n-1}{i} \\
& +&  \displaystyle{
    \sum_{i=0}^{n-1} }\bincc{m+n+\ell-k-1}{\ell-1}
        \bincc{\ell-1}{k-i-2}
    \bincc{m+n-i-2}{n-1}\bincc{n-1}{i}. 
\end{eqnarray*} }
On the other hand, using Pascal's identity
we see that the right hand side of the equation~(\ref{eq:3sh})
equals
\allowdisplaybreaks{
\begin{eqnarray*}
\lefteqn{ \displaystyle{
    \sum_{i=0}^n }\bincc{m+n+\ell-k}{\ell}\bincc{\ell}{k-i}
    \bincc{m+n-i}{n}\bincc{n}{i}=}\\
&& \displaystyle{
    \sum_{i=0}^n }\bincc{m+n+\ell-k-1}{\ell}\bincc{\ell}{k-i}
    \bincc{m+n-i}{n}\bincc{n}{i}
+ \displaystyle{
    \sum_{i=0}^n }\bincc{m+n+\ell-k-1}{\ell-1}\bincc{\ell}{k-i}
    \bincc{m+n-i}{n}\bincc{n}{i}\\
&=& \displaystyle{
    \sum_{i=0}^n }\bincc{m+n+\ell-k-1}{\ell}\bincc{\ell}{k-i}
    \bincc{m+n-i-1}{n}\bincc{n}{i}
+\displaystyle{
    \sum_{i=0}^n }\bincc{m+n+\ell-k-1}{\ell}\bincc{\ell}{k-i}
    \bincc{m+n-i-1}{n-1}\bincc{n}{i}\\
&&+ \displaystyle{
    \sum_{i=0}^n }\bincc{m+n+\ell-k-1}{\ell-1}
        \bincc{\ell-1}{k-i}
    \bincc{m+n-i-1}{n}\bincc{n}{i}
+ \displaystyle{
    \sum_{i=0}^n }\bincc{m+n+\ell-k-1}{\ell-1}
        \bincc{\ell-1}{k-i-1}
    \bincc{m+n-i}{n}\bincc{n}{i}\\
&=& \displaystyle{
    \sum_{i=0}^n }\bincc{m+n+\ell-k-1}{\ell}\bincc{\ell}{k-i}
    \bincc{m+n-i-1}{n}\bincc{n}{i}
+ \displaystyle{
    \sum_{i=0}^n }\bincc{m+n+\ell-k-1}{\ell}\bincc{\ell}{k-i}
    \bincc{m+n-i-1}{n-1}\bincc{n-1}{i}\\
&&+ \displaystyle{
    \sum_{i=0}^n }\bincc{m+n+\ell-k-1}{\ell}\bincc{\ell}{k-i}
    \bincc{m+n-i-1}{n-1}\bincc{n-1}{i-1}
+ \displaystyle{
    \sum_{i=0}^n }\bincc{m+n+\ell-k-1}{\ell-1}
        \bincc{\ell-1}{k-i}
    \bincc{m+n-i-1}{n}\bincc{n}{i}\\
&&+ \displaystyle{
    \sum_{i=0}^n }\bincc{m+n+\ell-k-1}{\ell-1}
        \bincc{\ell-1}{k-i-1}
    \bincc{m+n-i-1}{n}\bincc{n}{i}
+ \displaystyle{
    \sum_{i=0}^n }\bincc{m+n+\ell-k-1}{\ell-1}
        \bincc{\ell-1}{k-i-1}
    \bincc{m+n-i-1}{n-1}\bincc{n}{i}\\
&=& \displaystyle{
    \sum_{i=0}^n }\bincc{m+n+\ell-k-1}{\ell}\bincc{\ell}{k-i}
    \bincc{m+n-i-1}{n}\bincc{n}{i}
+ \displaystyle{
    \sum_{i=0}^n }\bincc{m+n+\ell-k-1}{\ell}\bincc{\ell}{k-i}
    \bincc{m+n-i-1}{n-1}\bincc{n-1}{i}\\
&&+ \displaystyle{
    \sum_{i=0}^n }\bincc{m+n+\ell-k-1}{\ell}\bincc{\ell}{k-i}
    \bincc{m+n-i-1}{n-1}\bincc{n-1}{i-1}
+  \displaystyle{
    \sum_{i=0}^n }\bincc{m+n+\ell-k-1}{\ell-1}
        \bincc{\ell-1}{k-i}
    \bincc{m+n-i-1}{n}\bincc{n}{i}\\
&&+  \displaystyle{
    \sum_{i=0}^n }\bincc{m+n+\ell-k-1}{\ell-1}
        \bincc{\ell-1}{k-i-1}
    \bincc{m+n-i-1}{n}\bincc{n}{i}
+ \displaystyle{
    \sum_{i=0}^n }\bincc{m+n+\ell-k-1}{\ell-1}
        \bincc{\ell-1}{k-i-1}
    \bincc{m+n-i-1}{n-1}\bincc{n-1}{i}\\
&&+ \displaystyle{
    \sum_{i=0}^n }\bincc{m+n+\ell-k-1}{\ell-1}
        \bincc{\ell-1}{k-i-1}
    \bincc{m+n-i-1}{n-1}\bincc{n-1}{i-1}\\
&=& \displaystyle{
    \sum_{i=0}^n }\bincc{m+n+\ell-k-1}{\ell}\bincc{\ell}{k-i}
    \bincc{m+n-i-1}{n}\bincc{n}{i}
+ \displaystyle{
    \sum_{i=0}^{n-1} }\bincc{m+n+\ell-k-1}{\ell}\bincc{\ell}{k-i}
    \bincc{m+n-i-1}{n-1}\bincc{n-1}{i}\\
&&+ \displaystyle{
    \sum_{i=0}^{n-1} }\bincc{m+n+\ell-k-1}{\ell}
        \bincc{\ell}{k-i-1}
    \bincc{m+n-i-2}{n-1}\bincc{n-1}{i}
+ \displaystyle{
    \sum_{i=0}^n }\bincc{m+n+\ell-k-1}{\ell-1}
        \bincc{\ell-1}{k-i}
    \bincc{m+n-i-1}{n}\bincc{n}{i}\\
&&+ \displaystyle{
    \sum_{i=0}^n }\bincc{m+n+\ell-k-1}{\ell-1}
        \bincc{\ell-1}{k-i-1}
    \bincc{m+n-i-1}{n}\bincc{n}{i}
+ \displaystyle{
    \sum_{i=0}^{n-1} }\bincc{m+n+\ell-k-1}{\ell-1}
        \bincc{\ell-1}{k-i-1}
    \bincc{m+n-i-1}{n-1}\bincc{n-1}{i}\\
&&+ \displaystyle{
    \sum_{i=0}^{n-1} }\bincc{m+n+\ell-k-1}{\ell-1}
        \bincc{\ell-1}{k-i-2}
    \bincc{m+n-i-2}{n-1}\bincc{n-1}{i}. 
\end{eqnarray*}}
This completes the induction, proving equation~(\ref{eq:3sh}). 
Then the equation in part 4 is obtained by summing
equation~(\ref{eq:3sh}) for $k=0,\ldots, n+\ell$. 
\proofend

\subsection{Mixable shuffles of vectors}

Let $\Omega$ be a countable infinite set. Let $\widetilde{\Omega}$ be the set
consisting of finite non-empty subsets of $\Omega$.
For $F=(F_1,\ldots,F_m)\in \widetilde{\Omega}^m$, $\sigma\in S_m$
and $T$ a subset of $\{ (\ell,\ell+1)\mid 1\leq l< m\}$,  
denote
\[ \sigma(F; T) = (F_{\sigma(1)}\hcm \ldots \hcm F_{\sigma(m)})\]
where, for each $1\leq \ell < m$, 
\[ F_{\sigma(\ell)}\hcm F_{\sigma(\ell+1)}=
\left \{ \begin{array}{ll} F_{\sigma(\ell)}\cup F_{\sigma(\ell+1)}
    & {\rm\ if\ } (\ell,\ell+1)\in T\\
    F_{\sigma(\ell)}, F_{\sigma(\ell+1)}
    & {\rm otherwise}
    \end{array}
\right . \]

For example, let $x_1,x_2,y,z\in \Omega$, $F_1=\{x_1,x_2\}$
$F_2=\{y\},\ F_3=\{z\}$,
$\sigma=\left ( \begin{array}{ccc} 1&2&3 \\
    2&1&3 \end{array} \right) \in S_3$,
$T=\{ (2,3) \}$. 
Then
\[ \sigma(F)= (F_{\sigma(1)}, F_{\sigma(2)}, F_{\sigma(3)})
    = (F_2,F_1,F_3)\]
and 
\allowdisplaybreaks{
\begin{eqnarray*}
\sigma(F; T) &=&
    (F_2 \hcm F_1 \hcm F_3)\\
    &=&
    (F_2 , F_1\cup F_3)\\
    &=&
    (\{y\}, \{x_1,x_2,z\}) \in \widetilde{\Omega}^2
\end{eqnarray*}}

Let $F=(F_1,\ldots,F_m)\in \widetilde{\Omega}^m,\ G=(G_1,\ldots,G_n)\in \widetilde{\Omega}^n$.
Denote $(F,G)=(F_1,\ldots,F_m,G_1,\ldots,G_n)\in \widetilde{\Omega}^{m+n}$
and denote
\[ U_k=\left \{ \begin{array}{ll}
    F_k, & 1\leq k\leq m\\
    G_{k-m}, & m+1\leq k\leq m+n
    \end{array} \right . \]

\begin{defn}
\begin{enumerate}
\item
For $\sigma\in S(m,n)$, $F\in \widetilde{\Omega}^m$, $G\in \widetilde{\Omega}^n$,
\[\sigma (F, G)=(U_{\sigma(1)},\ldots, U_{\sigma(m+n)})
    \in \widetilde{\Omega}^{m+n}\]
is called 
a {\bf shuffle} of $F$ and $G$.
\item
Denote 
\[ S(F,G)=\{ \sigma(F, G) \mid \sigma\in S(m,n)\}\]
for the {\bf set of shuffles} of $F$ and $G$. 
\item
Let $\sigma\in S(m,n)$ and
let $T$ be a subset of $\calt_\sigma$. 
The element 
\[ \sigma((F,G); T)= (U_{\sigma(1)}\hcm U_{\sigma(2)} \hcm
    \ldots \hcm U_{\sigma(m+n)}) \]
is called a {\bf mixable shuffle} of $F$ and $G$.
\item
Denote
\[ \bar{S}(F,G)=\{ \sigma((F, G); T) \mid
    \sigma\in S(m,n),\ T\subseteq \calt_\sigma \}. \]
for the {\bf set of mixable shuffles} of $F$ and
$G$. 
\end{enumerate}
\end{defn}

If we further have
$H=(H_1,\ldots,H_\ell)\in \widetilde{\Omega}^\ell$,
then we denote
\[ \bar{S}(\bar{S}(F,G),H)=\bigcup_{ U\in \bar{S}(F,G)}
     \bar{S}( U, H), \]
and     
\[ \bar{S}(F,\bar{S}(G,H))=\bigcup_{ U\in \bar{S}(G,H)}
     \bar{S}( F,U). \]
For $\sigma\in S(m,n,\ell)$ and $T\in \calt_\sigma$, 
define a mixable shuffle of $F,\ G$ and $H$ by 
\[ \sigma((F, G, H);T)=
    (U_{\sigma(1)}\hcm U_{\sigma(2)}\hcm \ldots \hcm
    U_{\sigma(m+n+\ell)}), \]
in  which 
\[ U_k=\left \{ \begin{array}{ll}
    F_k,& 1\leq k\leq m,\\
    G_{k-m}, & m+1\leq k\leq m+n \\
    H_{k-m-n},& m+n+1 \leq k\leq m+n+\ell. \end{array}
    \right .  \]
Then for $1\leq k< m+n+\ell-1$, we have 
\begin{eqnarray*}
\lefteqn{ U_{\sigma(k)}\hcm U_{\sigma(k+1)} \hcm U_{\sigma(k+2)}
    =}\\
&& \left \{\begin{array}{ll}
    U_{\sigma(k)}\cup U_{\sigma(k+1)}\cup U_{\sigma(k+2)}  &
     (k,k+1,k+2)\in T\\
    U_{\sigma(k)}\cup U_{\sigma(k+1)}, U_{\sigma(k+2)}  &
     (k,k+1,k+2)\not\in T, (k,k+1)\in T\\
    U_{\sigma(k)}, U_{\sigma(k+1)}\cup U_{\sigma(k+2)} &
     (k,k+1,k+2)\not\in T, (k+1,k+2)\in T\\
    U_{\sigma(k)}, U_{\sigma(k+1)}, U_{\sigma(k+2)}
    & {\rm otherwise}
    \end{array} \right .
\end{eqnarray*}
Denote 
\[\bs(F,G,H)= \{ \sigma((F, G, H); T) \mid
\sigma \in S(m,n,\ell),\ T\subseteq \calt_\sigma \}.\]
for the set of mixable shuffles of $F,\ G$ and $H$.

\begin{prop}
\label{prop:shuf}
Let $X_1,\ldots,X_m,Y_1,\ldots,Y_n,Z_1,\ldots,Z_\ell$ be
disjoint subsets of $\Omega$. 
Let $X=(X_i)\in \widetilde{\Omega}^m, \ Y=(Y_j)\in \widetilde{\Omega}^n,\
Z=(Z_k)\in \widetilde{\Omega}^\ell$.
Then 

\begin{enumerate}
\item
$\mid \bs(X,Y)\mid =s(m,n)$. 
\item
 $\bs(X,Y)=\bs(Y,X)$.
\item
$\bs(\bs(X,Y),Z)=\bs(X,Y,Z)=\bs(X,\bs(Y,Z)).$
\end{enumerate}
\end{prop}
\proof
1.
We use induction for $m+n$ with $m,n\geq 1$.
When $m=1$,
$S(X,Y)$ contains the vectors
\[ (X_1,Y_1,\ldots,Y_n), 
(Y_1,\ldots,Y_{i-1},X_1,Y_i, \ldots,Y_n), i=2,\ldots,n,
(Y_1,\ldots,Y_n,X_1)\]
and
\begin{eqnarray*}
&&\!\!\! (X_1\cup Y_1, Y_2,\ldots, Y_n),
(Y_1,\ldots,Y_{i-1},X_1\cup Y_i,Y_{i+1},\ldots,Y_n), i=2,\ldots,n-1, \\
&&  (Y_1,\ldots,Y_{n-1},X_1\cup Y_n).
\end{eqnarray*}
So there are at least $2n+1$ elements in $S(X,Y)$.
By construction,
the set $S(X,Y)$ has no more elements than
$S(m,n)$ which is $s(1,n)=2n+1$ by Proposition~\ref{prop:num}.
So the claim holds in this case.
A similar argument verifies the claim in the case when
$n=1$. 

Now assume that the claim is true for $m+n<k$ with $m,n> 1$,
and consider the case when $m+n=k$. 
Since $m,\ n\geq 2$, it makes sense to define 
\begin{eqnarray*}
 \bs_{X_1}(X,Y)&=& \{ (X_1 , S) \mid S\in
    \bs(X_2, \ldots X_m, Y) \}, \\
 \bs_{Y_1}(X,Y) &=&\{ (Y_1 , S) \mid
    S\in \bs(X,Y_2, \ldots , Y_n) \}, \\
 \bs_{X_1\cup Y_1} (X,Y) &=& \{  (X_1\cup Y_1 , S) \mid
    S\in \bs(X_2, \ldots , X_m,
        Y_2, \ldots , Y_n) \}.
\end{eqnarray*}
By assumption, $X_1,Y_1$ and $X_1\cup Y_1$ are distinct.
So we have 
\begin{equation}
 \bs(X,Y) \supseteq \bs_{X_1}(X,Y)\dotcup \bs_{Y_1}(X,Y) \dotcup
    \bs_{X_1\cup Y_1} (X,Y)
\label{eq:find}
\end{equation}
By induction hypothesis, the three sets on the right
have cardinalities $s(m-1,n), s(m,n-1)$ and $s(m-1,n-1)$.
So
\[ \mid \bs(X,Y)\mid\geq s(m-1,n)+s(m,n-1)+s(m-1,n-1)\]
which is $s(m,n)$ by equation~(\ref{eq:dec}).
So by Proposition~\ref{prop:num}, $\mid \bs(X,Y)\mid =s(m,n)$.
Then from equation~(\ref{eq:find}) we have
\begin{equation}
\bs(X,Y) = \bs_{X_1}(X,Y)\dotcup \bs_{Y_1}(X,Y) \dotcup
    \bs_{X_1\cup Y_1} (X,Y).
\label{eq:find2}
\end{equation}

2.
We prove by induction on $m+n,\ m,n\geq 1$. 
The statement can be directly verified for $m+n\leq 2$. 
Assume that it is true for $m+n<k$ and
let $X\in A^{m}$, $Y\in A^{n}$ with
$m+n=k$.
Then 
\begin{eqnarray*} \bs_{X_1}(X,Y)&=& \{ (X_1,S)\mid
    S\in \bs((X_2,\ldots,X_m),Y)\}\\
    &=& \{ (X_1,S)\mid
    S\in \bs(Y,(X_2,\ldots,X_m))\}\\
    &=& \bs_{X_1}(Y,X)
\end{eqnarray*}
and similarly,  
\[ \bs_{Y_1}(X,Y)=\bs_{Y_1}(Y,X)
{\rm\ and\ } \bs_{ X_1\cup Y_1}(X,Y)
=\bs_{Y_1\cup X_1}(Y,X). \]
Then it follows from Equation~(\ref{eq:find2}) and its symmetric form
for $\bs(Y,X)$ that $\bs(X,Y)=\bs(Y,X)$.

3. 
We first prove that
\[ \bs(\bs(X,Y),Z)=\bigcup _{U\in \bs(X,Y)} \bs(U,Z)\]
is a disjoint union.
By assumption, $X_i,Y_j$ and $Z_k$ are disjoint subsets of $\Omega$.
Since each component of $U=(U_r)\in \bs (X,Y)$ is either a $X_i$, or
$Y_j$ or a $X_i\cup Y_j$, it follows that the subset $U_r$
and $Z_k$ are also disjoint.
Let $U=(U_1,\ldots,U_r)$ and $U\tp=(U\tp_1,\ldots,U\tp_s)$
be two distinct mixable shuffles of $X$ and $Y$.
If $r\neq s$,
then without loss of generality, we could
assume that $r>s$.
Then there is a $U_{r_0}$ that is different from any $U\tp_j$
and therefore is disjoint with any $U\tp_j$.  
Thus $U_{r_0}$ is disjoint with any component of any
$W\in \bs(U\tp,Z)$.
On the other hand, $U_{r_0}$ has non-trivial intersection
with some component of every $W\in \bs(U,S)$.
Therefore, $\bs(U,Z)\cap \bs(U\tp,Z)=\phi$. 

Now assume that $r=s$.
We use induction on $r$.
For $r=1$, $U\neq U\tp$ means $U_1$ is different from
any components of $U\tp$.
So $U_1$ is disjoint from any components of $W\in \bs(U\tp,S)$, 
and the claim is proved.
Assume that the claim is true for $r$,
and let $U$ and $U\tp$ both have length $r+1$.
Suppose $\bs(U,Z)\cap\bs(U\tp,Z)$ is not empty.
Then there is a $W\in \bs(U,Z)\cap \bs(U\tp,Z)$. 
Write $W=(W_1,\cdots, W_k)$,
then $W_1=U_1$ or $Z_1$ or $U_1\cup Z_1$.
If $W_1=U_1$, then since $U_1$ is disjoint from any $Z_\ell$,
from $W\in \bs(U\tp,Z)$ we get $W_1=U\tp_1$.
This shows that

\[ (W_2,\ldots,W_k)\in \bs((U_2,\ldots,U_r),Z)
    \cap \bs((U\tp_2,\ldots,U\tp_r),Z).\]
But since $U_1=W_1=U\tp_1$, from $U\neq U\tp$ we get
$(U_2,\ldots,U_r)\neq (U\tp_2,\ldots,U\tp_r)$.
Then by induction assumption,

\[\bs((U_2,\ldots,U_r),Z)
    \cap \bs((U\tp_2,\ldots,U\tp_r),Z)=\phi.\]
This is a contradiction.
For the same reason, $W_1=Z_1$ or $U_1\cup Z_1$ also
implies contradiction.
Therefore, the claim is true for $r+1$.
This proves that
\[ \bs(\bs(X,Y),Z) =
\displaystyle{ \bigcup^\bullet_{U\in \bs(X,Y)}} \bs(U,Z).\]

By Proposition~\ref{prop:num}, $\bs(X,Y)$ has
$\bincc{n}{i} \bincc{m+n-i}{n}$ mixable shuffles of
length $i$, and for each such mixable shuffle $U$,
the cardinality of $\bs(U,Z)$ is
\[\sum_{j=0}^\ell \bincc{m+n+\ell-i-j}{\ell}
    \bincc{\ell}{j}.\]
Therefore $\bs (\bs(X,Y),Z)$ has
\begin{eqnarray*}
\lefteqn{\displaystyle{\sum_{i=0}^n
    \sum_{j=0}^\ell }
    \bincc{m+n-i}{n} \bincc{n}{i} 
    \bincc{m+n+\ell-i-j}{\ell} \bincc{\ell}{j}=}\\
&&\displaystyle{\sum_{k=0}^{n+\ell}
    \sum_{i=0}^n }\bincc{m+n+\ell-k}{\ell}\bincc{\ell}{k-i}
    \bincc{m+n-i}{n}\bincc{n}{i}
\end{eqnarray*}
elements.
But $\bs(\bs(X,Y),Z)$ is a subset of $\bs(X,Y,Z)$,
and
by Proposition~\ref{prop:num}, the second set has
cardinality bounded by the same number.
This proves the first equality in part 3 of the
proposition. 

To prove the second equality in part 3, 
define 
\[ \bs(X,Y)_{X_m} =\{ (S, X_m) \mid
    S\in \bs(X_1,\ldots , X_{m-1}, Y) \} \]
and similarly for $\bs(X,Y)_{Y_n}$,
$\bs (X,Y)_{X_m\cup Y_n}$.
Then the same argument as above gives 
\begin{equation}
\bs(X,Y)= \bs(X,Y)_{X_m}\dotcup \bs(X,Y)_{Y_n} \dotcup
    \bs(X,Y)_{X_1\cup Y_1}
\label{eq:eind}
\end{equation}
The rest of the proof is similar. 
\proofend

\subsection{Mixable shuffles of tensors}

\vspace{.3cm}
\noindent
{\bf Notation:} 
For the rest of this paper, let $\lambda$ be a fixed element of $C$. 
For any $C$-modules $M$ and $N$,
the tensor product $M\otimes N$ is taken over $C$
unless otherwise indicated. 
Let $A$ be a $C$-algebra. 
For $n\in \NN$, denote
\[ A^{\otimes n}=\underbrace{A\otimes\ldots\otimes A}
    _{n\ {\rm factors}}\]
with the convention that $A^{\otimes 0}=C$.

\vspace{.3cm}

Let $A$ be a $C$-algebra.
For $m,n,\ell\in \NN_+$, 
denote
$x=x_1\otimes\ldots\otimes x_m\in A^{\otimes m}$,
$y=y_1\otimes \ldots\otimes y_n\in A^{\otimes n}$
and
$z=z_1\otimes\ldots\otimes z_\ell\in A^{\otimes \ell}$.
For $\sigma\in S_n$, denote
\[\sigma (x)=
    x_{\sigma(1)}\otimes x_{\sigma(2)} \otimes
    \ldots \otimes x_{\sigma(n)}.\]
Denote
$x\otimes y =x_1\otimes \ldots  x_m\otimes y_1\ldots\otimes y_n
    \in A^{\otimes (m+n)}$. 
and, for $\sigma\in S_{m+n}$,
denote
\[\sigma (x\otimes y) =u_{\sigma(1)}\otimes u_{\sigma(2)} \otimes
    \ldots \otimes u_{\sigma(m+n)},\]
where
\[ u_k=\left \{ \begin{array}{ll}
    x_k,& 1\leq k\leq m,\\
    y_{k-m}, & m+1\leq k\leq m+n. \end{array}
    \right . \]
Likewise, for $\sigma\in S_{m+n+\ell}$,
denote
\[ \sigma(x\otimes y\otimes z)=
    u_{\sigma(1)}\otimes u_{\sigma(2)}\otimes \ldots \otimes
    u_{\sigma(m+n+\ell)}, \]
in  which 
\[ u_k=\left \{ \begin{array}{ll}
    x_k,& 1\leq k\leq m,\\
    y_{k-m}, & m+1\leq k\leq m+n \\
    z_{k-m-n},& m+n+1 \leq k\leq m+n+\ell. \end{array}
    \right .  \]

\begin{defn}
Let $x\in A^{\otimes m}$, $y\in A^{\otimes n}$
and $\sigma\in S(m,n)$.
\begin{enumerate}
\item
$\sigma (x\otimes y)\in A^{\otimes (m+n)}$ is called 
a {\bf shuffle} of $x$ and $y$.
\item
Let $T$ be a subset of $\calt_\sigma$. 
The element 
\[ \sigma(x\otimes y, T)= z_{\sigma(1)}\hts z_{\sigma(2)} \hts
    \ldots \hts z_{\sigma(m+n)}, \]
where for each pair $(k,k+1)$, $1\leq k< m+n$, 
\[ z_{\sigma(k)}\hts z_{\sigma(k+1)} =\left \{\begin{array}{ll}
    z_{\sigma(k)} z_{\sigma(k+1)},  &
     (k,k+1)\in T\\
    z_{\sigma(k)}\otimes z_{\sigma(k+1)}, &
    (k,k+1) \not\in T
    \end{array} \right . \]
is called a {\bf mixable shuffle} of $x$ and $y$.
\end{enumerate}
\end{defn}

It follows from the universal property of the tensor product
$A^{\otimes k}$ that $\sigma(x\otimes y; T)$ does not depend
on the choice of $x_1,\ldots,x_m,y_1,\ldots,y_n$
representing the tensor $x\otimes y$.

Now fix $\lambda\in C$. 
Define, for $x$ and $y$ as above, 

\[x \shpr\!\!^+ y\
=\sum_{(\sigma,T)\in \bs (m,n)} \lambda^{\mid T\mid } \sigma(x\otimes y;T) 
    \in \bigoplus_{k\leq m+n} A^{\otimes k}.\]
The operation $\shpr\!^+ $ extends to a mapping
\[ \shpr\!^+\ : A^{\otimes m}\times A^{\otimes n} \rar
    \bigoplus_{k\leq m+n} A^{\otimes k}, m, n\in \NN\]
by $C$-linearity.
Let

\[ \sha_C^+(A)=\sha_C^+(A,\lambda)= \bigoplus_{k\in\NN} A^{\otimes k}
=C\oplus A\oplus A^{\otimes 2}\oplus \ldots. \]
Extending by additivity, the binary operation
$\shpr^+$ gives a $C$-bilinear map 

\[ \shpr \!^+\ : \sha_C^+(A) \times \sha_C^+(A) \rar \sha_C^+(A) \]
with the convention that
\[ C\times A^{\otimes m} \rar A^{\otimes m} \]
is the scalar multiplication. 
This binary operation is called the
{\bf mixable shuffle product
of weight $\lambda$}.

\begin{theorem}
The mixable
shuffle product $\shpr \!^+$ defines 
an associative, commutative binary operation on
$\sha_C^+(A)= \bigoplus_{k\in\NN} A^{\otimes k}$, making it into
a $C$-algebra
with the identity $\bfone_C\in C=A^{\otimes 0}$.
\end{theorem}

$\sha_C^+(A)$ will be called the
{\bf mixable shuffle algebra (of weight $\lambda$)}
on $A$. 
In the special case when $\lambda=0$,
$\sha_C^+(A)$ is denoted by ${\rm Sh}(A)$ in~\cite{Sw}. 

\proof
We only need to verify the commutativity and associativity
of the operation $\shpr^+$.
For this we make use of the mixable shuffle of vectors discussed in
the previous section.

Recall that $\Omega$ is an infinite set,
and $\widetilde{\Omega}$ is the set
of finite non-empty subsets of $\Omega$.
For $F=(F_1,\ldots,F_m)\in \widetilde{\Omega}^m$, $\sigma\in S_m$
and $T$ a subset of $\{(\ell,\ell+1)\mid 1\leq \ell <m\}$, 
denote
\[ \sigma(F; T) = (F_{\sigma(1)}\hcm \ldots \hcm F_{\sigma(m)})\]
where, for each $1\leq \ell < m$, 
\[ F_{\sigma(\ell)}\hcm F_{\sigma(\ell+1)}=
\left \{ \begin{array}{ll} F_{\sigma(\ell)}\cup F_{\sigma(\ell+1)}
    & {\rm\ if\ } (\ell,\ell+1)\in T\\ 
    F_{\sigma(\ell)}, F_{\sigma(\ell+1)}   
    & {\rm otherwise}
    \end{array}
\right . \]

Given a $\varphi: \cup_{i=1}^m F_i \rar A$,
denote
\[ \varphi(\sigma(F; T)) =
 (a_{\sigma(1)}\hts \ldots \hts a_{\sigma(m)})\]
where, for each $1\leq \ell < m$, 
\[ a_{\sigma(\ell)}= \prod_{f_\ell\in F_\ell} \varphi(f_\ell).\]

For example, suppose $X_1,X_2,Y,Z\in \Omega$ are distinct.
Let $F_1=\{X_1,X_2\}$,
$F_2=\{Y\},\ F_3=\{Z\}$,
$\sigma=\left ( \begin{array}{ccc} 1& 2& 3 \\
    2&1&3 \end{array} \right) \in S_3$, $T=\{\ (2,3) \}$. 
If $\varphi$ sends each of the capital letters to the
corresponding lower case letter in the polynomial algebra
$A= C[x_1,x_2,y,z]$, then
\begin{eqnarray*}
\varphi(\sigma(F; T)) &=&
    \varphi (F_{\sigma(1)}, F_{\sigma(2)}\cup F_{\sigma(3)})\\
    &=&
    y\otimes x_1x_2 z \in A^{\otimes 2}. 
\end{eqnarray*}

For any fixed $m,n,\ell$ and fixed
$x\in A^{\otimes m},\ y\in A^{\otimes n}$ and
$z\in A^{\otimes \ell}$, choose distinct elements
$X_1,\ldots,X_m,Y_1,\ldots,Y_n,Z_1,\ldots,Z_\ell$ from
the infinite set $\Omega$. Also use the same letters for the
singletons $\{X_i\}, \{Y_j\}$ and $\{Z_k\}$ of $\Omega$.
Let $\varphi$ be the map sending $X_i$ to $x_i$, $Y_j$ to $y_j$
and $Z_k$ to $z_k$.
Then $\shpr^+$ (of weight $\lambda$) could be described as
 \[x \shpr\!\!^+ y = \varphi(X) \shpr^+ \varphi(Y) 
=\sum_{U\in \bs (X,Y)} \lambda^{\deg (U)} \varphi(U), \]
where $\deg (U)=\mid T\mid $ 
if $U$ is given by $(\sigma,T)\in \bs(m,n)$. 
The commutativity of $\shpr \!^+$ follows from
Proposition~\ref{prop:shuf}.2.

Let a mixable shuffle $W\in \bs(X,Y,Z)$ be given by
$(\sigma,T)\in \bs(m,n,\ell)$,
i.e., $W=\sigma((X,Y,Z);T)$.
Define $\deg (W)=\deg T $, with $\deg T$ defined in
equation~(\ref{eq:deg3}).
By Proposition~\ref{prop:shuf}.3,
$W$ could also be obtained from a mixable shuffle $V$ 
in $\bs(U,Z)$, given by $(\sigma_1,T_1)$,
where $U$ is from a mixable shuffle in $\bs(X,Y)$,
given by $(\sigma_2,T_2)$.
Thus we have
$V=\sigma_1((U,Z);T_1)$ and
$U=\sigma_2((X,Y); T_2)$.
It follows from the definition of $W$ and $\deg W$ that
the length of the vector $W$ is
$m+n+\ell- \deg (W)$.
Since $W$ is also given by $V=\sigma_1((U,Z);T_1)$, its length is also
given by
\[ {\rm\ length\ of\ } U + \ell - \deg (V)
= m+n- \deg(U) +\ell -\deg (V).\]
Thus we have $\deg (W)=\deg (U)+\deg (V)$. 
Then it follows from Proposition~\ref{prop:shuf}.3 that 
\allowdisplaybreaks{
\begin{eqnarray*}
(x\shpr^+ y)\shpr^+ z &=&
    (\sum_{U\in \bs(X,Y)} \lambda^{\deg (U)}\varphi(U))\shpr^+ z\\
    &=& \sum_{U\in \bs(X,Y)} \lambda^{\deg (U)}\varphi(U) \shpr^+ \varphi(Z)\\
    &=& \sum_{U\in \bs(X,Y)} \lambda^{\deg (U)} \sum_{V\in \bs(U,Z)}
         \lambda^{\deg (V)} \varphi(V)\\
    &=& \sum_{W \in \bs(\bs(X,Y),Z)}
        \lambda^{\deg (U)+\deg (V)}  \varphi(W)\\
    &=& \sum_{W\in \bs(X,Y,Z)} \lambda^{\deg (W)} \varphi(W). 
\end{eqnarray*}}
We similarly have
\[ x\shpr^+ (y\shpr^+ z)
    =\sum_{W\in \bs(X,Y,Z)} \lambda^{\deg (W)} \varphi(W). \]
This proves the associativity.
\proofend

\section{The free Baxter algebra}
\label{sec:free}

We now use the mixable shuffle product from last section
to describe the free Baxter algebras.
We will first construct the free Baxter algebra on
a $C$-algebra $A$.
We will then give constructions of other types of 
free Baxter algebras.

\subsection{The basic free construction}
With the same notations as those in last section,
we define $\sha_C(A)$ to be the tensor product algebra
$A\otimes_C \sha_C^+(A).$
Thus
\[ \sha_C(A)\cong \bigoplus_{k\in\NN} A^{\otimes (k+1)}
=A\oplus A^{\otimes 2}\oplus \ldots \]
as $A$-modules and 
the product on $\sha_C(A)$ is defined by the
{\bf augmented mixable shuffle product
(of weight $\lambda$)}
\begin{eqnarray*}
\lefteqn{(x_0\otimes x_1\otimes \ldots \otimes x_m)
\shpr (y_0\otimes y_1\otimes \ldots \otimes y_n)}\\
&=& x_0y_0 \otimes ((x_1\otimes \ldots\otimes x_m)\shpr\!\!^+
    (y_1\otimes \ldots\otimes y_n))\\
&=& x_0y_0 \otimes \sum_{(\sigma,T)\in \bs(m,n)}
\lambda^{\mid T\mid } \sigma ((x_1\otimes x_1\otimes \ldots \otimes x_m)
\otimes (y_1\otimes \ldots \otimes y_n); T)\\
&& \in A^{\otimes (m+n+1)}.
\end{eqnarray*}
Define a $C$-linear endomorphism $P_A$ on
$\sha_C(A)$ by assigning
\[ P_A( x_0\otimes x_1\otimes \ldots \otimes x_n)
=\bfone_A\otimes x_0\otimes x_1\otimes \ldots\otimes x_n, \]
for all
$x_0\otimes x_1\otimes \ldots\otimes x_n\in A^{\otimes (n+1)}$
and extending by additivity.
Let $j_A:A\rar \sha_C(A)$ be the canonical inclusion map.

\begin{theorem}
\label{thm:shua}
For any $C$-algebra $A$, 
$(\sha_C(A),P_A)$, together with the natural embedding
$j_A:A\rightarrow \sha_C(A)$,
is a free Baxter $C$-algebra on $A$ (of weight $\lambda$) in the sense that
the triple $(\sha_C(A),P_A,j_A)$ satisfies the following
universal property:
For any Baxter $C$-algebra $(R,P)$ and any
$C$-algebra map
$\varphi:A\rar R$, there exists
a unique Baxter $C$-algebra homomorphism
$\tilde{\varphi}:(\sha_C(A),P_A)\rar (R,P)$ such that
the  diagram
\[\xymatrix{
A \ar[rr]^(0.4){j_A} \ar[drr]_{\varphi}
    && \sha_C(A) \ar[d]^{\tilde{\varphi}} \\
&& R } \]
commutes. 
\end{theorem}

\remark
By the same argument used to show the uniqueness
of other ``universal'' objects, $\sha_C(A)$ is the unique
free Baxter $C$-algebra on $A$ up to isomorphism.
The Baxter $C$-algebra $\sha_C(A)$ will be
called the {\bf free Baxter $C$-algebra
(of weight $\lambda$) on $A$.}

\vspace{0.3cm}
\noindent
{\bf Proof: }
We first show that $P_A$ is a Baxter operator on $\sha_C(A)$. 
For this we only need to verify that for any $x,y\in \sha_C(A)$,
\[ P_A(x)\shpr P_A(y) = P_A(x\shpr P_A(y))
    +P_A(y\shpr P_A(x))+\lambda P_A(x\shpr y).\]
By additivity, we only need to verify this equation for
any 
$x=x_1\otimes \ldots\otimes x_m\in A^{\otimes m}$
and $y=y_1\otimes \ldots\otimes y_n\in A^{\otimes n}$.
By definition,
\begin{eqnarray*}
P_A(x)\shpr P_A(y) & =&
(\bfone_A\otimes x_1\otimes \ldots\otimes x_{m})\shpr
    (\bfone_A\otimes y_1\otimes \ldots\otimes y_{n})\\
&=& \bfone_A\otimes \sum_{ (\sigma,T)\in \bs(m,n)} \lambda^{\mid T\mid}
    \sigma(x\otimes y;T)\\
&=& P_A(\sum_{ (\sigma,T)\in \bs(m,n)} \lambda^{\mid T\mid}
    \sigma(x\otimes y;T)). 
\end{eqnarray*}

Recall that  equation~(\ref{eq:dec}) gives us 
\[
 \bs (m,n) = \bs_{1,0}(m,n) \dotcup
    \bs_{0,1}(m,n) \dotcup \bs_{1,1}(m,n)
\]
where 
\begin{eqnarray*}
 \bs_{1,0}(m,n)&=& \{(\sigma,T)\in \bs(m,n) \mid (1,2)\not\in T,
    \sigma^{-1}(1)=1 \},\\
 \bs_{0,1}(m,n)&=& \{(\sigma,T)\in \bs(m,n)\mid (1,2)\not\in T,
    \sigma^{-1}(m+1)=1 \},\\
 \bs_{1,1}(m,n)&=&\{ (\sigma,T)\in \bs(m,n)\mid (1,2)\in T\}.
\end{eqnarray*}
Further,
\begin{eqnarray*}
\sum_{(\sigma,T)\in \bs_{1,0}(m,n)} \lambda^{\mid T\mid} \sigma (x\otimes y)
    \!\!&=&\!\!
x_1\otimes\!\!\! \sum_{(\sigma,T)\in \bs(m-1,n)}
    \lambda^{\mid T\mid}
    \sigma ((x_2\otimes \ldots \otimes x_{m})\otimes y;T) \\
    &=&
x \shpr (\bfone_A\otimes y_1\otimes \ldots \otimes y_{n};T) \\
    &=&
x \shpr P_A(y).
\end{eqnarray*}
where for $m=1$,
$\sigma ((x_2\otimes \ldots \otimes x_{m})\otimes y;T)
=\sigma(y)$ (only $T=\phi$ is possible).
Similarly,
\[ \sum_{(\sigma,T)\in \bs_{0,1}(m,n)} \lambda^{\mid T\mid}
    \sigma (x\otimes y;T) =
P_A(x) \shpr y. \]
and
\begin{eqnarray*}
\lefteqn{ \sum_{(\sigma,T)\in \bs_{1,1}(m,n)}
    \lambda^{\mid T\mid }\sigma (x\otimes y;T)}\\
&=& \sum_{(\sigma,T)\in \bs(m-1,n-1)}
    \lambda x_1 y_1 \otimes \lambda^{\mid T\mid}
    \sigma((x_2\otimes \ldots \otimes x_{m-1})
    \otimes (y_2\otimes \ldots \otimes y_{n-1});T) \\
&=& \lambda P_A(x\shpr y).
\end{eqnarray*}
Therefore,

\[
P_A(x)\shpr P_A(y) =
P_A(x\shpr P_A(y)) + P_A(P_A(x)\shpr y)
+ \lambda P_A(x\shpr y)
\]
This shows that $P_A$ is a Baxter operator of weight $\lambda$ on
$\sha_C(A)$, making it into a Baxter $C$-algebra. 

Before verifying that $(\sha_C(A),P_A)$ satisfies
the universal property of a free Baxter $C$-algebra over $A$,
we need some preparations. 
Let $(R,P)$ be a Baxter $C$-algebra.
For $x,\ y\in R$, denote
$ P_x(y) = P(xy)$. 
For $x_0\otimes\ldots \otimes x_k\in R^{\otimes (k+1)},$
denote
\begin{equation}
 x_0 (\circ_{r=1}^k P_{x_r}) (y)=
    x_0 (P_{x_1}\circ\cdots \circ P_{x_k})(y),
\label{eq:comp}
\end{equation}
with the convention that
$\circ_{r=1}^0 P_{x_r}=\id_R.$
It follows from the universal property of the tensor product
$R^{\otimes (k+1)}$ 
and the $C$-linearity of the Baxter operator $P$ that
the right hand side of equation~(\ref{eq:comp})
is well-defined,
and does not depend on the choice of
$x_0,\ \ldots,\ x_k\in R$ representing the tensor
$x_0\otimes \ldots\otimes x_k$. 
For $\sigma\in S_n$, denote
\[\sigma (\circ_{r=1}^n P_{x_r})=
    \circ_{r=1}^n P_{x_{\sigma(r)}}.\]
For $m,\ n\in \NN_+$,
denote
\[(\circ_{r=1}^m P_{x_r}) \circ (\circ_{s=1}^n P_{y_s})
 = \circ_{t=1}^{m+n} P_{z_{t}},\]
where
\[ z_k=\left \{ \begin{array}{ll}
    x_k,& 1\leq k\leq m,\\
    y_{k-m}, & m+1\leq k\leq m+n \end{array}
    \right . \]
and, for $(\sigma,T)\in \bs(m,n)$,  denote
\[ \sigma((\circ_{r=1}^m P_{x_r} )\circ (\circ_{s=1}^n P_{y_s});T)
    =P_{z_1} \hcirc \ldots \hcirc P_{z_{m+n}}\]
where for each $(k,k+1), 1\leq k< m+n$, 
\[ P_{z_k}\hcirc P_{z_{k+1}} = \left \{
    \begin{array}{ll}
    P_{z_kz_{k+1}},& (k,k+1)\in T\\
    P_{z_k}\circ P_{z_{k+1}}, & (k,k+1)\not\in T
    \end{array}
\right .\]

\begin{prop}
\label{prop:prod1}
For $m, n\in \NN$ and
$x=x_0\otimes\ \ldots\otimes x_m\in R^{\otimes (m+1)},\
y=y_0\otimes \ldots\otimes y_n\in R^{\otimes (n+1)}$, 
\begin{eqnarray*}
\lefteqn{(x_0 (\circ_{r=1}^m P_{x_r})(\bfone_R))
    (y_0 (\circ_{s=1}^n P_{y_s})(\bfone_R))}\\
&& =x_0y_0\sum_{(\sigma,T)\in \bs(m,n)}
\lambda^{\mid T\mid}
\sigma ((\circ_{r=1}^m P_{x_r}) \circ (\circ_{s=1}^n P_{y_s}))
(\bfone_R).
\end{eqnarray*}
\end{prop}

\proof
It is clear that the equation in general follows from
the case when $x_0=y_0=1$, in which case the equation is 
\[(\circ_{r=1}^m P_{x_r}(\bfone_R))
    (\circ_{s=1}^n P_{y_s}(\bfone_R))
=\sum_{(\sigma,T)\in \bs(m,n)}
\lambda^{\mid T\mid}
\sigma ((\circ_{r=1}^m P_{x_r}) \circ (\circ_{s=1}^n P_{y_s}))
(\bfone_R)\]
with $m,\ n\geq 1.$ 

For this we prove by induction on $k=m+n$.
If $k=2$, then $m=n=1$. The equation to be proved in this case
is
\[ P(x_1)P(y_1)=P(x_1P(y_1))+P(y_1P(x_1))
    +\lambda P(x_1 y_1)\]
which is part of the definition of $P$. 
Assuming that the equation holds for all
$x\in R^{\otimes (m+1)},\ y\in R^{\otimes (n+1)}$ with $m+n<k$.
Let $x\in R^{\otimes (m+1)},\ y\in R^{\otimes (n+1)}$ with $m+n=k$.
Then
\allowdisplaybreaks{
\begin{eqnarray*}
&&  (\circ_{r=1}^m P_{x_r}) (\bfone_R)
     (\circ_{s=1}^n P_{y_s}) (\bfone_R)\\
&=& P(x_1 (\circ_{r=2}^m P_{x_r})(\bfone_R))
    P(y_1 (\circ_{s=2}^n P_{y_s})(\bfone_R))\\
&=&  P(x_1 (\circ_{r=2}^m P_{x_r})(\bfone_R)
    P(y_1 (\circ_{s=2}^n P_{y_s})(\bfone_R)))\\
&&  + P(y_1 (\circ_{s=2}^n P_{y_r})(\bfone_R)
    P(x_1 (\circ_{r=2}^m P_{x_r})(\bfone_R)))\\
&& + \lambda P(x_1 (\circ_{r=2}^m P_{x_r})(\bfone_R)
    y_1 (\circ_{s=2}^n P_{y_s})(\bfone_R))\\
&=&  P(x_1 (\circ_{r=2}^m P_{x_r})(\bfone_R)
    (\circ_{s=1}^n P_{y_s})(\bfone_R))
    + P((\circ_{r=1}^m P_{x_r})(\bfone_R)
    y_1 (\circ_{s=2}^{n} P_{y_r})(\bfone_R))\\
&& + \lambda P(x_1y_1 (\circ_{r=2}^m P_{x_r})(\bfone_R)
     (\circ_{s=2}^n P_{y_s})(\bfone_R))\\
&=& P(x_1 \sum_{(\sigma,T)\in \bs(m-1,n)}
    \lambda^{\mid T\mid}\sigma((\circ_{r=2}^m P_{x_r})\circ 
    (\circ_{s=1}^n P_{y_s}))(\bfone_R))\\
&& + P(y_1 \sum_{(\sigma,T)\in \bs(m,n-1)}
    \lambda^{\mid T\mid} \sigma((\circ_{r=1}^m P_{x_r})\circ 
    (\circ_{s=2}^n P_{y_s}))(\bfone_R))\\
&& +\lambda P(x_1y_1 \sum_{(\sigma,T)\in \bs(m-1,n-1)}
    \lambda^{\mid T\mid}\sigma((\circ_{r=2}^m P_{x_r})\circ 
    (\circ_{s=2}^n P_{y_s}))(\bfone_R))\\
&&  \ ({\rm by\ induction})\\
&=& \sum_{(\sigma,T)\in \bs(m-1,n)}
    \lambda^{\mid T\mid} P_{x_1}\circ  
    \sigma((\circ_{r=2}^m P_{x_r})\circ 
    (\circ_{s=1}^n P_{y_s}))(\bfone_R)\\
&& + \sum_{(\sigma,T)\in \bs(m,n-1)}
    \lambda^{\mid T\mid} P_{y_1}\circ  
    \sigma((\circ_{r=1}^m P_{x_r})\circ 
    (\circ_{s=2}^n P_{y_s}))(\bfone_R)\\
&& + \lambda \sum_{(\sigma,T)\in \bs(m-1,n-1)}
    \lambda^{\mid T\mid} P_{x_1y_1}\circ  
    \sigma((\circ_{r=2}^m P_{x_r})\circ 
    (\circ_{s=1}^n P_{y_s}))(\bfone_R)\\
&=& \sum_{(\sigma,T)\in \bs_{1,0}(m,n)}
    \lambda^{\mid T\mid}
    \sigma((\circ_{r=1}^m P_{x_r})\circ 
    (\circ_{s=1}^n P_{y_s}))(\bfone_R)\\
&& + \sum_{(\sigma,T)\in \bs_{0,1}(m,n)}
    \lambda^{\mid T\mid}
    \sigma((\circ_{r=1}^m P_{x_r})\circ 
    (\circ_{s=1}^n P_{y_s}))(\bfone_R)\\
&&+  \sum_{(\sigma,T)\in \bs_{1,1}(m,n)} 
    \lambda^{\mid T\mid}
    \sigma((\circ_{r=1}^m P_{x_r})\circ 
    (\circ_{s=1}^n P_{y_s}))(\bfone_R)\\
&=& \sum_{\sigma\in \bs(m,n)} 
    \lambda^{\mid T\mid }
    \sigma((\circ_{r=1}^m P_{x_r})\circ 
    (\circ_{s=1}^n P_{y_s}))(\bfone_R)\\
\end{eqnarray*}}
This completes the proof of Proposition~\ref{prop:prod1}.
\proofend

We now continue with the proof of Theorem~\ref{thm:shua}
and verify the universal property for $\sha_C(A)$. 
For a given Baxter $C$-algebra $A$ and
a $C$-algebra map $\varphi:A\rar R$, 
we extend $\varphi$ to an Baxter $C$-algebra homomorphism
$\tilde{\varphi}: \sha_C(A)\rar (R,i)$ as follows. 
For any 
$x=x_0\otimes x_1\otimes \ldots\otimes x_k\in A^{\otimes (k+1)}$,
define
\[\tilde{\varphi}(x)
= \varphi(x_0) (\circ_{j=1}^k P_{\varphi(x_j)}) (\bfone_R) .\]
This is a well-defined $C$-linear map, hence extends uniquely by
additivity to a $C$-module homomorphism
\[\tilde{\varphi}: \sha_C(A)= \bigoplus_{k\in\NN} A^{\otimes (k+1)} 
    \rar R  .\]
It follows from the definition of the operation $\shpr$ and
Proposition~\ref{prop:prod1} that $\tilde{\varphi}$
preserves multiplication. 
Since $\tilde{\varphi}$ is $C$-linear, and 
\begin{eqnarray*}
\tilde{\varphi}(P_A(x_0\otimes \ldots\otimes x_k))
&=& \tilde{\varphi}(\bfone_A\otimes x_0\otimes \ldots\otimes x_k))\\
&=& (\circ_{j=0}^k P_{\varphi(x_j)})(\bfone_R)\\
&=& P(\varphi(x_0) (\circ_{j=1}^k P_{\varphi(x_j)})(\bfone_R)) \\
&=& P(\tilde{\varphi}(x_0\otimes \ldots\otimes x_k)), 
\end{eqnarray*}
$\tilde{\varphi}$ is a homomorphism of Baxter $C$-algebras.
It is clear from its construction that
$\tilde{\varphi}$ is the unique
homomorphism of Baxter $C$-algebras extending $\varphi$.
This proves Theorem~\ref{thm:shua}.  
\proofend

Let $\Alg_C$ be the category of $C$-algebras and
let $U_C:\Bax_C\rar \Alg_C$ be the forgetful
functor.  

\begin{coro}
\label{co:adj}
\begin{enumerate}
\item
The assignment $A\mapsto \sha_C(A)$ gives a functor 
\[ \sha_C: \Alg_C\rar \Bax_C,\]
where for a $C$-algebra homomorphism $f:A\rar B$,
$\sha_C(f): \sha_C(A)\rar \sha_C(B)$ is defined by
$\oplus_{n\in \NN} f^{\otimes n}$ with
$f^{\otimes n}:A^{\otimes n}\rar B^{\otimes n}$. 
\item
$\sha_C$ is the 
left adjoint functor of the forgetful functor
$U_C$. 
\item
Any Baxter $C$-algebra is isomorphic to a quotient of
$(\sha_C(A),P_A)$ for some $C$-algebra $A$.
\end{enumerate}
\end{coro}

By Corollary~\ref{co:adj}, the study of Baxter $C$-algebras 
is reduced to studying quotients of
free Baxter $C$-algebras.

{ } From the proof of the theorem and Proposition~\ref{prop:prod1},
we also obtain

\begin{coro}
\label{co:gen}
\begin{enumerate}
\item
For any sub-$C$-algebra $B$ of a Baxter $C$-algebra $(R,P)$,
the Baxter sub-$C$-algebra of $R$ generated by $B$ is
generated by
\[ \{ x_0 (\circ_{i=1}^r P_{x_i})(\bfone_R)\mid x_i\in B,\
    0\leq i\leq r,\ r\in \NN \}\]
as an additive group.
\item
For any $C$-algebra homomorphism $\varphi:A\rar R$,
the image of $\tilde{\varphi}:\sha_C(A)\rar R$
is the Baxter subalgebra of $(R,P)$ generated by
$\varphi(A)$. 
\end{enumerate}
\end{coro}

\proof
Let $\hat{B}$ be the Baxter sub-$C$-algebra of $(R,P)$
generated by $B$. 
Denote 
\[ S= \{ x_0 (\circ_{i=1}^r P_{x_i})(\bfone_R)\mid x_i\in B,\
    0\leq i\leq r,\ r\in \NN \}\]
and let $\hat{B}\tp$ be the additive subgroup of $R$ generated by $S$.
Since $S$ is closed under scalar multiplication by $C$ and
the operator $P$, 
$\hat{B}\tp$ is a $C$-module and is closed under $P$. 
By Proposition~\ref{prop:prod1}, $\hat{B}\tp$ is closed under
multiplication. Therefore, $\hat{B}\tp$ is a sub-$C$-algebra of $R$,
hence contains $\hat{B}$.
On the other hand, since $\hat{B}$ contains $B$ and
is closed under multiplication and Baxter operator,  it must
contain $S$. Then $\hat{B}$ must contain $\hat{B}\tp$
by closure under addition.
This proves the first statement.

By its construction, the Baxter $C$-algebra $\sha_C(A)$
is generated by $j_A(A)$.
Since $\tilde{\varphi}$ is an Baxter $C$-algebra
homomorphism, $\tilde{\varphi}(\sha_C(A))$ is
also an Baxter $C$-subalgebra, and is generated by
$\tilde{\varphi}(j_A(A))=\varphi(A)$.
\proofend

\subsection{Other free constructions}
\label{sec:var}
The construction of the free Baxter $C$-algebra
$\sha_C(A)$ 
on a $C$-algebra $A$ in the last part
could be combined with other free constructions to obtain
free Baxter algebras on other structures.
We now discuss
the free Baxter algebra on a set. 
The free Baxter algebras on a monoid or on a $C$-module
will also be considered.

For a given set $X$, 
let $C[X]$ be the polynomial $C$-algebra on $X$
with the natural embedding $X\hookrightarrow C[X]$. 
Let $(\sha_C(X),P_X)$ be the Baxter $C$-algebra
$(\sha_C(C[X]),P_{C[X]})$.

\begin{prop}
\label{co:shux}
$(\sha_C(X),P_X)$, together with the set embedding
\[j_X: X\hookrightarrow C[X] \ola{j_{C[X]}} \sha_C(C[X]),\]
is a free Baxter $C$-algebra on the set $X$,
described by the following universal property:
For any Baxter $C$-algebra $(R,P)$ over $C$ and any set map
$\varphi: X\rar R$, there exists
a unique Baxter $C$-algebra homomorphism
$\tilde{\varphi}:(\sha_C(X),P_X)\rar (R,i)$ such that
the diagram 
\[\xymatrix{
X \ar[rr]^(0.4){j_X} \ar[drr]_{\varphi}
    && \sha_C(X) \ar[d]^{\tilde{\varphi}} \\
&& R } \]
commutes. 
\end{prop}

\remark
When $\lambda=-1$, it is easy to show that
$\sha_C(X)$ is closely related to the free Baxter
algebra constructed by Cartier\cite{Ca}.
See section~\ref{sec:RC} for detail.

\vspace{0.3cm}
\noindent
\proof
\noindent
We only need to verify that $(\intg{X},P_X)$ satisfies
the universal property of a free Baxter $C$-algebra
on $X$. 
Fix a given Baxter algebra $(R,P)$ over $C$, and a set 
map $\varphi:X\rar R$.
By the universal property of the $C$-algebra
$C[X]$, $\varphi$ extends uniquely by multiplicity
and $C$-linearity to a
$C$-algebra homomorphism
$\bar{\varphi}:C[X] \rar R$.
Then by Theorem~\ref{thm:shua}, $\bar{\varphi}$
extends uniquely to a homomorphism of Baxter $C$-algebras
\[\tilde{\varphi}: (\sha_C(X),P_X)\cong
    (\sha_C(C[X]),P_{C[X]}) \rar (R,i).\]
This proves the proposition. 
\proofend

Note that the free Baxter $C$-algebra on a set $X$
can be described as the composite of two free constructions.
First we construct the free $C$-algebra $C[X]$ on $X$, and then
we construct the free Baxter $C$-algebra on $C[X]$.
In a similar manner, we can construct the free Baxter
$C$-algebra on a monoid $M$ by first constructing the
free $C$-algebra $C\modg{M}$ on $M$
and then the
free Baxter $C$-algebra on $C\modg{M}$.
Here, for a given commutative monoid $M$, 
$C\modg{M}$ is the free $C$-module $\oplus_{x\in M} Cx$, 
where the multiplication on $C\modg{M}$ is induced
by the multiplication on $M$.

Likewise, we can construct the free Baxter
$C$-algebra on a $C$-module $N$ by first constructing the
tensor $C$-algebra $T_C(N)$ on $N$ and then the
free Baxter $C$-algebra on $T_C(N)$.

\section{Relation to constructions of Cartier and Rota}
\label{sec:RC}

We now give an explicit description of the relation between
our construction of
free Baxter algebras and the constructions of Cartier~\cite{Ca}
and Rota~\cite{Rot1}.

As before, let $C$ be a commutative ring with identity. 
Let $\Alg^0_C$ be the category in which the objects are
$C$-algebras not necessarily having an identity and the morphisms
preserve the addition and the multiplication, but do not
necessarily preserve the identity. 
Rota and Cartier considered the category 
$\Bax^0_C$ whose objects
are pairs $(R,P)$ where $R$ is an object in $\Alg_C^0$ and $P$ is
a Baxter operator (of weight $-1$), and whose morphisms
are morphisms in $\Alg^0_C$ that commute with the Baxter operators.

For any non-empty set $X$, the existence of a free Baxter algebra on
$X$ in either of the two categories $\Bax_C$ or $\Bax_C^0$
can be proved by general results from universal algebra
\cite{Co,Ma}. 
Rota and Cartier have given explicit descriptions of the
free Baxter algebra in $\Bax_C^0$.
Below we will focus on the relation of our construction of
free Baxter algebras in $\Bax_C$ 
with Cartier's construction of free Baxter algebras 
in $\Bax_C^0$. 
This will also explain the relation with 
Rota's construction of free Baxter algebras in $\Bax_C^0$, 
since, by the uniqueness of universal objects in $\Bax_C^0$,
the free Baxter algebras of Cartier are 
isomorphic to the free Baxter algebras of Rota.

We first recall the free Baxter algebra in $\Bax_C^0$
constructed by Cartier.
Let $M$ be the free commutative semigroup with identity on $X$.
Let $\widetilde{X}$ denote the set of symbols of the form
\begin{eqnarray*}
 u_0\cdot [\ ] &,&\ u_0\in M,\ u_0\neq 1, {\rm\ and\ }\\
 u_0 \cdot [u_1,\ldots,u_m] &,& m\geq 1,\ u_0, u_1,\ldots,u_m\in M,\
    u_m\neq 1.
\end{eqnarray*}
Let $\frakB(X)$ be the free $C$-module on $\widetilde{X}$. 
Cartier gave a $C$-bilinear multiplication
$\shprc$ 
on $\frakB(X)$ by defining
\begin{eqnarray*}
(u_0 \cdot [\ ])\shprc (v_0\cdot [\ ])&=& u_0v_0 \cdot [\ ],\\
(u_0 \cdot [\ ])\shprc (v_0 \cdot [v_1,\ldots,v_n])
    &=& (v_0\cdot [v_1,\ldots,v_n])\shprc (u_0\cdot [\ ])\\
    &=& u_0 v_0\cdot [v_1,\ldots,v_n],
\end{eqnarray*}
and
\begin{eqnarray*}
\lefteqn{
(u_0 \cdot [u_1,\ldots,u_m])\shprc (v_0 \cdot [v_1,\ldots,v_n])}\\
    &=& \sum_{(k,P,Q)\in \bar{S}_c(m,n)}
    (-1)^{k+p+q} u_0 v_0 \cdot
    \Phi_{k,P,Q}([u_1,\ldots,u_m],[v_1\ldots,v_n]).
\end{eqnarray*}
Here $\bar{S}_c(m,n)$ is the set of triples $(k,P,Q)$ 
in which $k$ is an integer between $1$ and $m+n$, 
$P$ and $Q$ are ordered subsets of $\{1,\ldots, k\}$ with the
natural ordering such that
$P\cup Q=\{1,\ldots,k\}$, $\mid\!\! P\!\!\mid\ =m$ and
$\mid\!\! Q\!\!\mid\ =n$. 
For each $(k,P,Q)\in \bar{S}_c(m,n)$, 
$\Phi_{k,P,Q}([u_1,\ldots,u_m],[v_1,\ldots,v_n])$
is the element $[c_1,\ldots,c_k]$ in $\widetilde{X}$ defined by

\[
c_j=\left \{ \begin{array}{ll}
    a_\alpha, & {\rm\ if\ } j {\rm\ is\ the\ }
    \alpha-{\rm th\ element\ in\ } P,\ j\not\in Q;\\
    b_\beta, & {\rm\ if\ } j {\rm\ is\ the\ }
    \beta-{\rm th\ element\ in\ } Q,\ j\not\in P;\\
    a_\alpha b_\beta, & {\rm\ if\ } j {\rm\ is\ the\ }
    \alpha-{\rm th\ element\ in\ } P    \\
    & {\rm\ and\ }
    {\rm\ the\ } \beta-{\rm th\ element\ in\ } Q
    \end{array} \right . \]
Define a $C$-linear operator $P_X^c$ on $\frakB(X)$ by 
\begin{eqnarray*}
P_X^c(u_0\cdot[\ ])&=&1\cdot [u_0], \\
P_X^c(u_0\cdot [u_1,\ldots,u_m])&=&1\cdot [u_0,u_1,\ldots,u_m].
\end{eqnarray*} 
Cartier proved that the pair $(\frakB(X),P_X^c)$ is a
free Baxter algebra on $X$ in the category $\Bax_C^0$
    
On the other hand, since $C[X]$ is a free $C$-module on
$M$, it follows from our construction of the mixable shuffle
product algebra $\sha_C(X)=\sha_C(C[X])$ that
$\sha_C(X)$ is a free $C$-module on the set $\hat{X}$
of tensors
\[ u_0\otimes \ldots u_k,\ u_k\in M, k\geq 0. \]
The Baxter operator $P_X$ on $\sha_C(X)$ is defined by

\[P_X(u_0\otimes\ldots\otimes u_m)
=1\otimes u_0\otimes \ldots\otimes u_m. \]

We define a map $f: \widetilde{X}\to \hat{X}$ by
\begin{eqnarray*}
f(u_0\cdot [\ ])& =& u_0;\\
f(u_0\cdot [u_1,\ldots,u_m])
&=& u_0\otimes u_1\otimes\ldots\otimes u_m, 
\end{eqnarray*}
and extend it by $C$-linearity to a $C$-linear map

\[ f: \frakB(X)\to \sha_C(X). \]

\begin{prop}
\label{prop:inj}
$f$ is an injective morphism in $\Bax^0_C$,
identifying $\frakB(X)$ with the sub-Baxter algebra
of $\sha_C(X)$ with the $C$-basis
\[ \{ u_0\otimes\ldots u_m\mid m\geq 0,\ u_i\in M,\ i=0,\ldots, m,\
    u_m\neq 1\}.\]
\end{prop}

\proof
Define $h: X\to \frakB(X)$ by $h(x)=x\cdot [\ ]$.    
Cartier showed that
the pair $(\frakB(X),h)$ is a free Baxter algebra on $X$ in
the category $\Bax^0_C$.
Define $g: X\to \sha_C(X)$ by $g(x) =x\in \sha_C(X),\ x\in X$,
and regard $\sha_C(X)$ as an element in $\Bax^0_C(X)$. 
By the universal property of $\frakB(X)$,
there is a unique morphism
\[ \tilde{g}: \frakB(X)\to \sha_C(X) \]
in $\Bax_C^0$ such that
$\tilde{g}(x\cdot [\ ])= x$.
It follows from the fact that $\tilde{g}$ preserves multiplication
that
$\tilde{g}(u_0\cdot [\ ])=u_0$ for $u_0\in M,\ u_0\neq 1$.

It remains to prove,
for $m\geq 1,\ u_0,\ldots,u_m\in M,\ u_m\neq 1$, 
that 
\begin{equation}
 \tilde{g}(u_0\cdot [u_1,\ldots,u_m])
    = u_0\otimes u_1\otimes\ldots\otimes u_m.
\label{eq:inj}
\end{equation}
We will prove this by induction. 
Let $m=1$.  If $u_0=1$,  we have

\[ \tilde{g}(u_0\cdot [u_1]) =\tilde{g}(P_X^c(u_1))
=P_X(\tilde{g}(u_1))=1\otimes u_1.\]
If $u_0\neq 1$, we have
\begin{eqnarray*}
 \tilde{g}(u_0\cdot [u_1])
&=& \tilde{g}((u_0\cdot [\ ]) \shprc (1\cdot [u_1])) \\
&=& \tilde{g}(u_0\cdot [\ ])\shpr \tilde{g}(1\cdot [u_1])\\
&=& \tilde{g}(u_0\cdot [\ ])\shpr \tilde{g}(P_X^c(u_1\cdot [\ ]))\\
&=& u_0 \shpr  P_X(\tilde{g}(u_1\cdot [\ ]))\\
&=& u_0 \shpr P_X(u_1)\\
&=& u_0 \shpr (1\otimes u_1)\\
&=& u_0\otimes u_1.
\end{eqnarray*}
Assume that equation~(\ref{eq:inj}) holds for $m\geq 1$ and consider
the element 
\[ u_0\cdot [u_1,\ldots,u_{m+1}],\ u_i\in M, u_{m+1}\neq 1. \]
If $u_0 =1$, then
\begin{eqnarray*}
 \tilde{g}(u_0\cdot [u_1,\ldots,u_{m+1}])
&=&\tilde{g}(P_X^c(u_1\cdot [u_2,\ldots,u_{m+1}]))\\
&=& P_X(\tilde{g}(u_1\cdot [u_2,\ldots,u_{m+1}]))\\
&=& P_X(u_1\otimes \ldots \otimes u_{m+1})\\
&=& 1\otimes u_1\otimes \ldots \otimes u_{m+1}\\
&=& u_0\otimes u_1\otimes \ldots \otimes u_{m+1}
\end{eqnarray*}
If $u_0\neq 1$, we have
\begin{eqnarray*}
 \tilde{g}(u_0\cdot [u_1,\ldots,u_{m+1}])
&=& \tilde{g}((u_0\cdot [\ ]) \shprc (1\cdot [u_1,\ldots,u_{m+1}]))\\
&=& \tilde{g}(u_0\cdot [\ ])\shpr
    \tilde{g}(P_X^c(u_1\cdot [u_2,\ldots,u_{m+1}]))\\
&=& u_0 \shpr  P_X(\tilde{g}(u_1\cdot [u_2,\ldots,u_{m+1}]))\\
&=& u_0 \shpr P_X(u_1\otimes u_2\otimes\ldots\otimes u_{m+1})\\
&=& u_0 \shpr (1\otimes u_1\otimes\ldots\otimes u_{m+1})\\
&=& u_0 \otimes u_1\otimes\ldots\otimes u_{m+1}.
\end{eqnarray*}
Thus $\tilde{g}$ is the $C$-linear map $f$ defined above,
proving that $f$ is a morphism in $\Bax_C^0$.
Since $f$ is $C$-linear and
sends the $C$-basis $\tilde{X}$ of $\frakB(X)$
injectively to the $C$-linearly independent set

\[\{ u_0\otimes \ldots\otimes u_m\mid
    m\geq 0,\ u_i\in M, u_m\neq 1\},\]
we see that $f$ is injective.
\proofend

Proposition~\ref{prop:inj} enables us to identify
$\frakB(X)$ as a sub-Baxter algebra of $\sha_C(X)$
in the category $\Bax_C^0$. We further have 

\begin{prop}
\label{prop:bsu}
The injective morphism $f:\frakB(X)\to \sha_C(X)$ in
$\Bax_C^0$ satisfies
the following universal property.
For any element $A$ in $\Bax_C$, also regarded as an
element in $\Bax_C^0$, and any morphism
$\phi: \frakB(X) \to A$ in $\Bax_C^0$,
there is a unique morphism
$\tilde{\phi}: \sha_C(X)\to A$ in $\Bax_C$ such that
$\tilde{\phi}\circ f =\phi$.
\end{prop}

\proof
With the notations introduced in the proof of
Proposition~\ref{prop:inj}, 
we have the following diagram

\[\begin{array}{ccc}
 X &\ola{h} & \frakB(X) \\
 \dap{g} & {}^f\!\!\swarrow & \dap{\phi}\\
 \sha_C(X) & \ola{\tilde{\phi}} & A
\end{array}
\]
We only need to find $\tilde{\phi}: \sha_C(X)\to A$
in $\Bax_C$ such that the lower right triangle commutes.
By the universal property of $\sha_C(X)$ in $\Bax_C$,
the map $\phi\circ h: X\to A$ induces a morphism
$\psi: \sha_C(X)\to A$ in $\Bax_C$ such that
$\psi\circ g = \phi \circ h$.
We then have

\[ \psi\circ f\circ h = \psi\circ g = \phi\circ h.\]
From the universal property of $\frakB(X)$ in
$\Bax_C^0$, we obtain 
\[ \psi \circ f = \phi.\]
So we can take $\tilde{\phi}=\psi$.
\proofend

\section{The free Baxter $C$-algebra on $C$}
\label{sec:C}
As a particular example, we consider $\sha_C(C)$,
the free Baxter $C$-algebra on $C$.
$\sha_C(C)$ is also $\sha_C(\phi)$,
the free Baxter $C$-algebra on the empty set $\phi$,
defined in Corollary~\ref{co:shux}. 
This free Baxter algebra not only provides the simplest example
of a free Baxter algebra, it also helps to explain an
important difference
between our construction of the free Baxter algebra and
the construction of the free Baxter algebra of Rota or Cartier.
As we see from the previous section,
the free Baxter algebras of Rota and Cartier are
in $\Bax_C^0$ and have no identity.
In fact, the free Baxter algebra on $\phi$ in
$\Bax_C^0$ is the zero algebra with the zero Baxter operator.

If we choose $A=C$ in the construction of free Baxter
$C$-algebras, then we get
\[ \sha_C(C)=\bigoplus_{n=0}^\infty C^{\otimes (n+1)}. \]
Since the tensor product is over $C$, we have
$C^{\otimes (n+1)} =C \bfone^{\otimes (n+1)}$, 
where
$\bfone^{\otimes (n+1)}
= \underbrace{\bfone_C \otimes \ldots \otimes \bfone_C}
_{(n+1)-{\rm factors}}$.
In particular, $\bfone^{\otimes 1} = \bfone_C$. 
Thus $\sha_C(C)$ is a free $C$-module on the basis
$\bfone^{\otimes n}, n\geq 1$.

\begin{prop}
\label{prop:unit}
For any $m,n\in \NN$,
\[ \bfone^{\otimes (m+1)} \shpr \bfone^{\otimes (n+1)} =
\sum_{k=0}^m \binc{m+n-k}{n}\binc{n}{k} \lambda^k
\bfone^{\otimes (m+n+1-k)}.\]
\end{prop}

\proof
This is immediate from part 2 of
Proposition~\ref{prop:num}.
\proofend

We also have the following consequence of Theorem~\ref{thm:shua}.
\begin{coro}
\label{coro:init}
$\sha_C(C)$ is the initial object in the category $\Bax_C$ of
Baxter $C$-algebras. In other words, for any
Baxter $C$-algebra $(R,P)$, there is a unique
Baxter $C$-algebra homomorphism
$(\sha_C(C),P_C)\rar (R,P)$. 
\end{coro}

Now we consider the special case when $\lambda=0$. 
In this case this $C$-algebra has been studied earlier~\cite[\S 12.3]{Sw},
as the shuffle algebra $Sh(C)$ over $C$.
The following facts from there can easily be verified. 

\begin{prop}
\label{prop:sw}
Let $\lambda=0$. 
\begin{enumerate}
\item
The free Baxter $C$-algebra $\sha_C(C)$ is a bialgebra.
\item
For any $m,\ n\in \NN$, 
$ \bfone^{\otimes (m+1)} \shpr \bfone^{\otimes (n+1)}
= \binc{m+n}{n} \bfone^{\otimes (m+n+1)}.$
\end{enumerate}
\end{prop}

Now let $HC$ be the ring of Hurwitz series over $C$~\cite{Ke},
defined to be the set of sequences 
\[ \{ (a_n) \mid a_n\in C, n\in \NN \} \]
in which the addition
is defined componentwise and the multiplication is
defined by
\[ (a_n) (b_n) =(c_n)\]
with
\[ c_n =\displaystyle{ \sum_{k=0}^n }
    \bincc{n}{k}a_kb_{n-k}.\]
Denote $e_n$ for the sequence $(a_k)$ in which
$a_n=\bfone_C$ and $a_k=0$ for $k\neq n$.
Since $e_n e_m = \bincc{m+n}{n}e_{m+n}$,
from Proposition~\ref{prop:sw}, we obtain

\begin{prop}
Let $\lambda=0$. 
The assignment
\[ \bfone^{\otimes (n+1)} \mapsto e_n,\ n\geq 0\]
defines an injective homomorphism of
Baxter $C$-algebras from $\sha_C(C)$ to
$HC$, identifying $\sha_C(C)$ with the subalgebra
of ``Hurwitz polynomials''
\[ 
\{ (a_n)\in HC \mid a_n= 0, n>\!\!>0\}
= \{ \sum_{n\geq 0} a_n e_n \mid a_n =0, n>\!\!> 0\}. 
\]
\end{prop}

\addcontentsline{toc}{section}{\numberline {}References}

\end{document}